\documentclass[11pt]{amsart}   	
\usepackage[margin = 1in]{geometry}
									
\usepackage{amsmath,amssymb, amsfonts,amsthm}
\usepackage[dvipsnames]{xcolor}
\usepackage{hyperref}
\usepackage{mathtools}
\usepackage{tikz}
\usepackage{tikz-cd}
\usetikzlibrary{decorations.markings}
\usepackage{caption}
\usepackage{array}

\newtheorem{theorem}{Theorem}[section]
\newtheorem{lemma}[theorem]{Lemma}
\newtheorem{corollary}[theorem]{Corollary}
\newtheorem{proposition}[theorem]{Proposition}

\theoremstyle{definition}
\newtheorem{definition}[theorem]{Definition}
\newtheorem{example}[theorem]{Example}
\newtheorem{remark}[theorem]{Remark}

\DeclareMathOperator{\brick}{\mathsf{brick}}
\DeclareMathOperator{\sbrick}{\mathsf{sbrick}}
\newcommand{\arc}{\mathrm{arc}}
\DeclareMathOperator{\mods}{\mathsf{mod}}
\newcommand{\Hom}{\mathrm{Hom}}
\newcommand{\Ext}{\mathrm{Ext}}
\newcommand{\weak}{\mathrm{weak}}

\newcommand{\wrd}{\mathrm{wd}}

\newcommand{\undim}{\underline{\mathrm{dim}}}

\usepackage{soul}

\author{Eric J. Hanson}
\address{LACIM, Universit\'e du Qu\'ebec \`a Montr\'eal and Universit\'e de Sherbrooke, Qu\'ebec, CANADA}
\email{ejhanso3@ncsu.edu}

\author{Xinrui You}
\address{School of Mathematical Sciences, Nankai University, Tianjin, CHINA}
\email{yxr2001@icloud.com}

\subjclass[2020]{16G20, 16G10\\\indent This is the accepted manuscript of a paper published in {\it Journal of Algebra}. The official version of record is available at \url{https://doi.org/10.1016/j.jalgebra.2023.10.019}.  This manuscript version is available under a Creative Commons CC BY-NC-ND 4.0 license \url{https://creativecommons.org/licenses/by-nc-nd/4.0}. \copyright \ 2024.}

\title[Bricks in preprojective type A]{Morphisms and extensions between bricks over preprojective algebras of type A}

\date{8 November 2023}

\begin{document}
\maketitle

\begin{abstract}
    The bricks over preprojective algebras of type A are known to be in bijection with certain combinatorial objects called ``arcs''. In this paper, we show how one can use arcs to compute bases for the Hom-spaces and first extension spaces between bricks. We then use this description to classify the ``weak exceptional sequences'' over these algebras. Finally, we explain how our result relates to a similar combinatorial model for the exceptional sequences over hereditary algebras of type A.
\end{abstract}

\setcounter{tocdepth}{1}
\tableofcontents

\section{Introduction}

Let $Q$ be a quiver, $K$ a field, and $\Lambda = KQ/I$ the quotient of the path algebra $KQ$ by an admissible ideal $I$. An \emph{exceptional sequence} over $\Lambda$ is a sequence $(X_k,\ldots,X_1)$ of finitely generated (right) $\Lambda$-modules such that $\Hom_\Lambda(X_i,X_j) = 0$ for $1 \leq i < j \leq k$, $\Hom_\Lambda(X_i,X_i)$ is a division algebra (that is, $X_i$ is a \emph{brick}) for $1 \leq i \leq k$, and $\Ext^m_\Lambda(X_i,X_j) = 0$ for $1 \leq i \leq j \leq k$ and $m \in \mathbb{N}$. It is well-known that the length of an exceptional sequence
cannot exceed the number of vertices of $Q$. Furthermore, if $Q$ is acyclic and $I = 0$ (that is, if $\Lambda$ is hereditary) then every exceptional sequence can be ``completed'' so as to achieve this length, see \cite{CB_exceptional,ringel_exceptional}.

On the other hand, there are many examples of algebras admitting exceptional sequences which cannot be completed. For example, in this paper we focus on \emph{preprojective algebras} of Dynkin type A (see Definition~\ref{def:preproj}), denoted $\Pi(A_n)$. Over these algebras, we have that $\Hom_{\Pi(A_n)}(X,Y) \neq 0$ for any projective modules $X$ and $Y$ and that $\Ext^2_{\Pi(A_n)}(X,X) \neq 0$ for any non-projective finitely-generated module $X$, see e.g. \cite[Proposition~5.3]{reiten}. As a consequence, no exceptional sequence over $\Pi(A_n)$ can have length more than one. For additional information on the historical development and significance of preprojective algebras, we refer to the introduction of \cite{GLS} and the references therein.

Since the exceptional sequences over an arbitrary algebra may not have very nice properties, one natural approach is to weaken the definition. Indeed, to recover the definition in the hereditary case, one needs only assume that a sequence $(X_k,\ldots,X_1)$ of indecomposable modules satisfies $\Hom_\Lambda(X_i,X_j) = 0$ for $1 \leq i < j \leq k$ and $\Ext^1_\Lambda(X_i,X_j) = 0$ for $1 \leq i \leq j \leq k$. Over arbitrary algebras, sequences which satisfy only these properties are known as \emph{stratifying systems}. We refer to the introduction of \cite{MT} and the references therein for information about the rich history and body of research related to stratifying systems.

The preprojective algebra $\Pi(A_n)$ has wild representation type for $n \geq 6$. Thus classifying all stratifying systems over such algebras is a difficult question. On the other hand, the number of bricks over $\Pi(A_n)$ remains finite for any choice of $n$. This follows from the ``brick-$\tau$-rigid correspondence'' of \cite{DIJ} together with Mizuno's result that the ``poset of $\tau$-tilting pairs"of a preprojective algebra of Dynkin type is isomorphic to the weak order on the corresponding Coxeter group \cite{mizuno}. Thus it is natural to consider only those stratifying systems consisting entirely as bricks. Such stratifying systems are also known as \emph{weak exceptional sequences} \cite{sen}.

The main result of this paper (Theorem~\ref{thm:weak}) uses ``arcs'' on $n+1$ nodes to formulate a concrete combinatorial description of the weak exceptional sequences over the preprojective algebra $\Pi(A_n)$ (and certain quotient algebras). As an application, we use this model to prove that $2n-2$ is a tight upper bound on the length of a weak exceptional sequence over $\Pi(A_n)$ for $n \geq 2$ (Corollary~\ref{cor:planar2}). We also explain how our model relates to that of the exceptional sequences over path algebras of type $A_n$ established in \cite{GIMO} (Section~\ref{sec:hereditary}).

In concurrent work by the first author \cite{hanson}, we also use arcs to classify which weak exceptional sequences over $\Pi(A_n)$ are the sequences of bricks corresponding (in the sense of \cite[Theorem~D]{BaH2}) to Buan and Marsh's \emph{$\tau$-exceptional sequences} \cite{BM_exceptional}. These are another example of stratifying systems which, in full generality, satisfy a weaker version of the ``completion property'' observed for classical exceptional sequences in the hereditary case. See \cite{BM_exceptional}, \cite[Theorem~2.10]{AIR}, and \cite[Section~5]{MT} for further details.

Our results fit into a larger body of work which uses arcs (and associated combinatorial objects) to model various sets of representation-theoretic objects over $\Pi(A_n)$. This includes the sets of $\tau$-tilting modules \cite{IRRT,mizuno}, bricks and semibricks \cite{asai,BCZ,enomoto}, and 2-term simple minded collections \cite{BaH,mizuno2}. Arcs are also used to construct an explicit basis for the Hom-space between bricks in \cite{BCZ,mizuno2}. In constructing our model, we reprove this result (Proposition~\ref{prop:hom_space}) and also use arcs to give an explicit basis for the first Ext-space between two bricks (Theorem~\ref{thm:ext_space}). For a certain quotient of the preprojective algebra, this basis serves as a parallel of the description for when $\Hom(X,\tau Y)$ vanishes formulated in \cite[Section~3.2]{IW}. It also allows us to construct a model of the ``almost rigid modules'' (in the sense of \cite{BCSGS,BGMS}) over the same quotient, see Remark~\ref{rem:MAR}.

\subsection*{Funding and Acknowledgements}

X.Y. was supported by a Mitacs Globalink Summer Research Internship held at l'Universit\'e du Qu\'ebec \`a Montr\'eal. E.H. was supported by the Canada Research Chairs program (CRC-2021-00120) and NSERC Discovery Grants (RGPIN-2022-03960 and RGPIN/04465-2019). The authors are thankful to Emily Barnard, Karin Baur, Aslak Bakke Buan, Raquel Coehlo Sim{\~o}es, Benjamin Dequ{\^e}ne, Hern\'an Ibarra-Mejia, Ray Maresca, Clara Otte, and Hugh Thomas for many insightful discussions related to this project.


\section{Bricks and weak exceptional sequences}\label{sec:background}

Let $\Lambda$ be a finite-dimensional basic algebra over a field $K$. We denote by $\mods\Lambda$ the category of finitely-generated right $\Lambda$-modules. We adopt the common convention of identifying $\mods\Lambda$ with a skeleton, meaning we identify every $M \in \mods\Lambda$ with its isomorphism class.

A module $X \in \mods\Lambda$ is called a \emph{brick} if $\mathrm{End}_\Lambda(X)$ is a division algebra. When $K$ is algebraically closed, this condition simplifies to $\mathrm{End}_\Lambda(X) \cong K$. (This will also be the case for the type A preprojective algebras considered in the paper, see Corollary~\ref{cor:bricks_K}.) A set of bricks $\mathcal{X}$ is called a \emph{semibrick} if $\Hom_\Lambda(X,Y) = 0 = \Hom_\Lambda(Y,X)$ for all $X \neq Y \in \mathcal{X}$. We denote by $\brick(\Lambda)$ and $\sbrick(\Lambda)$ the sets of bricks and semibricks in $\mods\Lambda$, respectively.

The following is the main object of study in this paper.

\begin{definition}\cite{sen}\label{def:weak}
    Let $\Delta = (X_k,\ldots,X_1)$ be a sequence of indecomposable modules. We say that $\Delta$ is a \emph{weak exceptional sequence} if the following hold.
    \begin{enumerate}
        \item $\mathrm{End}_\Lambda(X_i) \cong K$ (so in particular $X_i$ is a brick) for all $1 \leq i \leq k$.
        \item $\Hom_\Lambda(X_i,X_j) = 0$ for all $1 \leq i < j \leq k$.
        \item $\Ext^1_\Lambda(X_i,X_j) = 0$ for all $1 \leq i \leq j \leq k$.
    \end{enumerate}
    We refer to $k$ as the \emph{length} of $\Delta$. In case $k = 2$, we will sometimes refer to the weak exceptional sequence $(X_2,X_1)$ as a \emph{weak exceptional pair}.
\end{definition}

\begin{remark}\label{rem:weak}\
    \begin{enumerate}
        \item If one replaces (3) in Definition~\ref{def:weak} with (3') $\Ext^m_\Lambda(X_i,X_j) = 0$ for all $1 \leq i \leq j \leq k$ and for all $m \geq 1$, the result is the definition of an \emph{exceptional sequence}. Similarly, if one removes condition (1) from Definition~\ref{def:weak}, the result is the definition of a \emph{stratifying system}. Moreover, if $\Lambda$ is hereditary then these three notions all coincide. As mentioned in the introduction, we refer to the introduction of \cite{MT} and the references therein for further discussion on these alternatives. 
        \item Weak exceptional pairs are an example of the ``semibrick pairs'' introduced in \cite{HI}. Indeed, if every brick in $\mods\Lambda$ has trivial endomorphism ring, then a pair $(X_2,X_1)$ is a weak exceptional sequence if and only if $X_1 \sqcup X_2[1]$ is a semibrick pair in the sense of \cite[Definition~1.8]{HI}. Over the type-A preprojective algebras studied in this paper, such pairs were classified in \cite[Section~6.2]{BaH} (see also Lemma~\ref{lem:same_weak} in the present paper).
    \end{enumerate}
\end{remark}

Note that a sequence $(X_k,\ldots,X_1)$ is a weak exceptional sequence if and only if $(X_j,X_i)$ is a weak exceptional pair for all $1 \leq i < j \leq k$. We denote by $\weak(\Lambda, k)$ the set of weak exceptional sequences in $\mods\Lambda$ of length $k$ and by $\weak(\Lambda) = \bigcup_{k \in \mathbb{N}} \weak(\Lambda,k)$. Note also that, outside of the hereditary case, the length of a weak exceptional sequence (and thus also of a stratifying system) may surpass the number of simple objects in $\mods\Lambda$. See e.g. \cite[3.2]{ES}, \cite[Remark~2.7]{MMS}, \cite{sen,treffinger}, or Corollary~\ref{cor:planar2} for examples.


\section{Preprojective algebras and their quotients}\label{sec:preproj}

 In this section, we give background information on preprojective algebras of type $A_n$. We begin with the definition. Our convention is to compose arrows from left-to-right.

\begin{definition}\label{def:preproj}
    Let $A_n$ be a linear quiver of type $A_n$, with an arrow $a_i: i \rightarrow i+1$ for $1 \leq i < n$. We denote by $\overline{A_n}$ the \emph{double quiver} of $A_n$, obtained by adjoining an arrow $a_i^*: i+1 \rightarrow i$ for each arrow $a_i$ in $Q$. We then define a family of algebras as follows. Let $S \subseteq [2,n-1]_\mathbb{Z} := \{2,3,\ldots,n-1\}$. We denote
        \begin{eqnarray*}
            c_2(S) &:=& \{a_1a_1^*,a_{n-1}^*a_{n-1}\} \cup \left(\bigcup_{i \in S}\{a_{i}a_{i}^*, a_{i-1}^*a_{i-1}\}\right) \cup \left(\bigcup_{i \in [2,n-1]_\mathbb{Z}\setminus S} \{a_{i}a_{i}^*-a_{i-1}^*a_{i-1}\}\right),\\
            \Pi(A_n,S) &:=& \overline{A_n}/(c_2(S)).
        \end{eqnarray*}
    Then $\Pi(A_n) := \Pi(A_n,\emptyset)$ is the \emph{preprojective algebra} of type $A_n$.
\end{definition}

\begin{example}
Consider the quivers 
\[\begin{tikzcd}[sep=scriptsize]
	A_4: & 1\arrow[r,"a_1"] & 2\arrow[r,"a_2"] & 3\arrow[r,"a_3"] & 4, && \overline{A_4}: & 1\arrow[r,"a_1",shift left] & 2\arrow[r,"a_2",shift left]\arrow[l,"a_1^*",shift left] & 3\arrow[r,"a_3",shift left]\arrow[l,"a_2^*",shift left] & 4.\arrow[l,"a_3^*",shift left]
\end{tikzcd}\]
Then the preprojective algebra $\Pi(A_4)$ and the algebra $\Pi(A_n,\{2,3\})$ are
\begin{eqnarray*}
\Pi (A_4) = \Pi(A_n,\emptyset) &=& K\overline{A_4}/ (a_1a_1^*,a_2a_2^*-a_1^*a_1, a_3a_3^*-a_2^*a_2,a_3^*a_3),\\
\Pi(A_4,\{2,3\}) &=& K\overline{A_4}/ (a_1a_1^*,a_2a_2^*,a_1^*a_1, a_3a_3^*,a_2^*a_2,a_3^*a_3).
\end{eqnarray*}
\end{example}

\begin{remark}
    While we have defined the preprojective algebra only for one particular quiver, the definition is in fact much more general. Indeed, for an arbitrary quiver $Q = (Q_0,Q_1)$, the \emph{preprojective algebra of type $Q$}, denoted $\Pi(Q)$, is defined to be the quotient of $K\overline{Q}$ by the ideal generated by $\sum_{a \in Q_1} a_ia_i^*-a_i^*a_i$. This algebra will be finite-dimensional if and only if $Q$ is of Dynkin type. Furthermore, it is well-known that (up to isomorphism) the algebra $\Pi(Q)$ depends only on the underlying (undirected) graph of $Q$, and not on the orientations of the arrows. Thus our choice to work only with the linear orientation of $A_n$ serves only to fix which arrows are denoted $a_i$ and which are denoted $a_i^*$.
\end{remark}

 Note that for $S \subseteq T \subseteq [2,n-1]_\mathbb{Z}$, one has $(c_2(S))\supseteq (c_2(T))$, and so there is a natural quotient map $q_T^S: \Pi(A_n,S) \rightarrow \Pi(A_n,T)$ which induces a fully faithful functor $F_T^S: \mods \Pi(A_n,T) \rightarrow \mods \Pi(A_n,S)$. When viewing modules as representations of quivers with relations, the functor $F_T^S$ acts as the identity on both objects and morphisms. In particular, for $S \subseteq T \subseteq U \subseteq [2,n-1]_\mathbb{Z}$ we have $F_U^S = F_U^T \circ F_T^S$. Moreover, each functor $F_T^S$ is exact and reflects exact sequences; that is, for $M, N \in \mods \Pi(A_n,S)$, the functor $F_T^S$ induces an injective linear map $\Ext^1_{\Pi(A_n,T)}(M,N) \hookrightarrow \Ext^1_{\Pi(A_n,S)}(M,N)$. This will be useful in Section~\ref{sec:ext} for proving that the elements of our proposed basis for the Ext-space between bricks in $\mods\Pi(A_n)$ are linearly independent.

The algebra $\Pi(A_n)$ is known to have wild representation type for $n \geq 6$. At the same time, it was shown by Mizuno \cite{mizuno} that $\Pi(A_n)$ is ``$\tau$-tilting finite'', and so in particular $\mods\Pi(A_n)$ contains only finitely many bricks by \cite{DIJ}. On the other hand, the algebra $\Pi(A_n,[2,n-1]_\mathbb{Z})$ is a ``gentle algebra'' with no ``bands''. Thus the classical work of Butler and Ringel \cite{butler_ringel} implies that $\mods \Pi(A_n,[2,n-1]_\mathbb{Z})$ has finitely many indecomposables, and that each of these is a brick.

In \cite{mizuno}, Mizuno showed that the ``lattice of torsion classes'' of $\Pi(A_n)$ is isomorphic to the ``weak order'' on the corresponding Coxeter group. In \cite{IRRT}, Iyama, Reading, Reiten, and Thomas use this to model the ``$\tau^{-1}$-rigid modules'' in $\mods\Pi(A_n)$ using ``Young-like diagrams''. This was modified to a model of the bricks over $\Pi(A_n)$ by Asai \cite{asai} via the dual of the brick-$\tau$-rigid correspondence of \cite{DIJ}. In the same year, Barnard, Carroll, and Zhu \cite[Section~4]{BCZ} and Demonet, Iyama, Reading, Reiten, and Thomas \cite[Section~6.3]{DIRRT} independently studied the bricks over $\Pi(A_n,[2,n-1]_\mathbb{Z})$. They showed that each brick can be identified with an orientation of a connected subgraph of $A_n$, and vice versa. Combinatorially, this phenomenon is modeled in \cite{BCZ} using ``arcs'', which had appeared in Reading's previous work on ``join-irreducible elements'' of the weak order \cite{reading}. The use of arcs to model bricks over $\Pi(A_n)$ is also established in \cite{enomoto,mizuno2} without explicit reference to the algebra $\Pi(A_n,[2,n-1]_\mathbb{Z})$. Furthermore, it is shown in both \cite{BCZ} and \cite{DIRRT} that the lattice of torsion classes of $\Pi(A_n,[2,n-1]_\mathbb{Z})$ is also isomorphic to the weak order on the corresponding Coxeter group (see also \cite{kase}), a consequence of which is the following.

\begin{proposition}\label{prop:sameBricks}
    Let $S \subseteq T \subseteq [2,n-1]_\mathbb{Z}$. Then $F_T^S$ induces an equality $\brick(\Pi(A_n,T)) = \brick(\Pi(A_n,S))$.
\end{proposition}

In the remainder of this section, we recall the ``arc model'' of $\brick(\Pi(A_n))$, which we identify with $\brick(\Pi(A_n),S)$ for all $S$ by the above proposition. We begin by recalling the definition of an arc on $n+1$ nodes.

\begin{definition}\label{def:arc}
    Let $[n] := \{0,\ldots,n\}$ be a set of $n+1$ nodes in $\mathbb{R}^2$, arranged in increasing order at the points $(i,0)$. An \emph{arc} $\gamma$ on $n+1$ nodes is the data of a pair of nodes $l(\gamma) < r(\gamma)$ and a continuous function $[l(\gamma),r(\gamma)] \rightarrow \mathbb{R}$ (which by abuse of notation we also denote by $\gamma$) such that
    \begin{enumerate}
        \item $\gamma(l(\gamma)) = 0 = \gamma(r(\gamma))$, and
        \item $\gamma(k) \neq 0$ for all $k \in (l(\gamma),r(\gamma)) \cap \mathbb{Z}$.
    \end{enumerate}
    We refer to $l(\gamma)$ and $r(\gamma)$ as the left and right \emph{endpoints} of $\gamma$, respectively.
\end{definition}

Alternatively, one can view an arc as a path $\gamma: [0,1] \rightarrow \mathbb{R}^2$ which moves monotonically from left to right and takes values in $[n]$ only at $t = 0$ and $t = 1$.

\begin{definition}\label{def:support}
    Let $\gamma$ be an arc on $n+1$ nodes. The \emph{support} of $\gamma$ is the open interval $(l(\gamma),r(\gamma)) \subseteq \mathbb{R}$ and the \emph{closed support} of $\gamma$ is the closed interval $[l(\gamma),r(\gamma)] \subseteq \mathbb{R}$. Likewise the \emph{arrow support} of $\gamma$ is the set $[l(\gamma) + 1,r(\gamma)-1]_\mathbb{Z} \subseteq [n]$. For each $k$ in the arrow support of $\gamma$, we say that $\gamma$ passes \emph{above} $k$ if $\gamma(k) > 0$ and that $\gamma$ passes \emph{below} $k$ if $\gamma(k) < 0$. Given a second arc $\rho$, the \emph{common (closed, arrow) support} of $\gamma$ and $\rho$ is the intersection of their (closed, arrow) supports.
\end{definition}

The name arrow support is justified in Remark~\ref{rem:visualize}

When considering an arc, we are mostly interested in its endpoints and whether it passes above or below each node $v$ in its arrow support. We encode this information as follows.

\begin{definition}\label{def:wd}
    Let $\mathrm{word}(\text{ueo})$ be the set of nonempty words on the letters u, e, and o. We treat $\{\text{u}, \text{e}, \text{o}\}$ as a totally ordered set with $\text{u} < \text{e} < \text{o}$. Define a map $\wrd:\arc(n) \rightarrow \{0,\ldots,n-1\} \times \mathrm{word}(\text{uoe})$ as follows. Let $\gamma \in \arc(n)$. For $l(\gamma) \leq i \leq r(\gamma)$, set
    $$s_i := \begin{cases} \text{u} & \text{if }\gamma(i) < 0\\ \text{o} & \text{if }\gamma(i) > 0\\ \text{e} & \text{if } \gamma(i) = 0.\end{cases}$$
    We then set $\wrd(\gamma) := \left(l(\gamma),s_{l(\gamma)+1}\cdots s_{r(\gamma)}\right).$
\end{definition}

\begin{figure}
\begin{tikzcd}[sep=scriptsize]  
        & K & K & K & K
	\arrow["1" below, shift left=1, from=1-3, to=1-2]
	\arrow["0" above, shift left=1, dashed, from=1-2, to=1-3]
	\arrow["1", shift left=1, from=1-3, to=1-4]
	\arrow["0", shift left=1, dashed, from=1-4, to=1-3]
	\arrow["1" below, shift left=1, from=1-5, to=1-4]
	\arrow["0" above, shift left=1, dashed, from=1-4, to=1-5]
\end{tikzcd}
\hspace*{2em}
$\vcenter{\hbox{
\begin{tikzpicture}[scale=0.8]
\draw [very thick] 
(-4,0) .. controls (-3.25, 1.3)and(-1.75,1.3) .. (-1.1,0)

(-1.1,0) .. controls (-0.5,-1.2)and(0.5, -1.2) .. (1.1,0)

          .. controls (1.75, 1.3)and(3.25,1.3) .. (4,0)
          node[near end, above]{};
\filldraw  (0,0) circle (2.4pt)
 (-2,0) circle (2.4pt)
 (-4,0) circle (2.4pt)
 (2,0) circle (2.4pt)
 (4,0) circle (2.4pt) 
 ;

\draw [very thick, ->] 
(-0.3, 0.22) .. controls (-0.15, 0.5)and(0.15, 0.5) .. (0.3, 0.22);

\draw [very thick, <-] 
(2-0.3, -0.22) .. controls (2-0.15, -0.5)and(2.15, -0.5) .. (2.3, -0.22);

\draw [very thick, <-]
(-2-0.3, -0.22) .. controls (-2-0.15, -0.5)and(0.15-2, -0.5) .. (0.3-2, -0.22);
\end{tikzpicture}}}$
\caption{A brick over $\Pi(A_n)$ and the corresponding arc on $n+1$ nodes.}\label{fig:arc}
\end{figure}
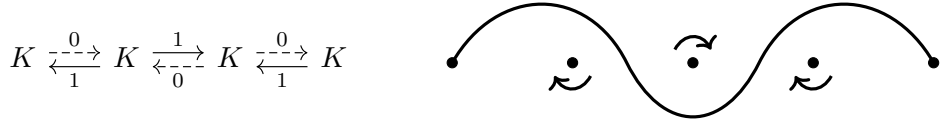

Definition~\ref{def:wd} can be used to rigorously say that two arcs $\gamma$ and $\rho$ are \emph{combinatorially equivalent} if and only if $\wrd(\gamma) = \wrd(\rho)$. We then denote by $\arc(n)$ the set of equivalence classes of arcs on $n+1$ nodes. We will regularly adopt the standard convention of using $\gamma \in \arc(n)$ to mean both the arc $\gamma$ and its equivalence class. To aid in visualization and help clarify the relationship between arcs and bricks over $\Pi(A_n)$, we sometimes decorate an arc $\gamma$ with \emph{small clockwise arrows} as shown in Figure~\ref{fig:arc}.

We now recall the bijection between bricks and arcs.

\begin{proposition}\cite[Proposition~4.6]{BCZ}\cite[Proposition~3.7]{mizuno2}\label{prop:arcBricks} Let $S \subseteq [2,n-1]_\mathbb{Z}$. Then:
    \begin{enumerate}
        \item There is a bijection $\sigma: \brick(\Pi(A_n,S))\rightarrow \arc(n)$ given as follows. Let $X \in \brick(\Pi(A_n,S))$. Then $X$ satisfies the following.
        \begin{enumerate}
            \item There exists $i \leq j \in [n]$ such that $\dim(X(k)) = 1$ for $i \leq k \leq j$ and $\dim(X(k)) = 0$ for $k < i$ or $k > j$.
            \item If $\dim(X(i)) = 1 = \dim(X(i+1))$, then one of $X(a_i)$ and $X(a_i^*)$ is 0 and the other is an isomorphism.
        \end{enumerate}
        We denote an arc $\sigma(X)$ as follows. Let $i, j$ be as in (a). Then $\wrd(\sigma(X)) = \left(i-1,s_i\cdots s_j\right)$ with
        $$s_k = \begin{cases} \textnormal{o} & \text{if } X(a_k^*) \neq 0\\\textnormal{u} & \text{if } X(a_k) \neq 0 \\ \textnormal{e} & \text{if }k = j.\end{cases}$$ That is, the endpoints of $\sigma(X)$ are $l(\sigma(X)) = i-1$ and $r(\sigma(X)) = j$. For each $i < k < j$, $\sigma(X)$ passes above the node $k$ if $X(a_k) \neq 0$ and passes below the node $k$ otherwise.

        \item Let $\gamma \in \arc(n)$ and denote $\wrd(\gamma) = \left(l(\gamma),s_{l(\gamma)+1}\cdots s_r\right)$. Then the brick $\sigma^{-1}(\gamma)$ is given by $$\sigma^{-1}(\gamma)(k) = \begin{cases} K & l(\gamma) < k \leq r(\gamma)\\0 & \text{ otherwise,}\end{cases}$$$$\sigma^{-1}(\gamma)(a_k) = \begin{cases} 1_K & s_k = \textnormal{u}\\0 & \text{otherwise},\end{cases}\qquad\qquad \sigma^{-1}(\gamma)(a_k^*) = \begin{cases} 1_K & s_k = \textnormal{o}\\0 & \text{otherwise}.\end{cases}$$
    \end{enumerate}
\end{proposition}

\begin{remark}\label{rem:visualize}
    Visually, we can understand the bijection $\sigma$ as follows. See Figure~\ref{fig:arc} for an example using the brick $X = \text{\tiny $\begin{matrix} \ \ 2 \ 4\\ 1 \ 3 \ \end{matrix}$} \in \brick(\Pi(A_4,S))$. The middle $n-1$ nodes represent the (pairs of) arrows of the quiver $\overline{A}_n$. Given a brick $X \in \brick(\Pi(A_n,S))$, we draw a small clockwise arrow to represent the direction of the nonzero map near its corresponding nodes. Then we draw an arc $\sigma(X)$ from node $i$ to $j$, where all the nonzero linear maps lie exactly between $i$ and $j$. When we have a small clockwise arrow above (resp.\ below) a node, the arc should go below (resp.\ above) the node. Note in particular that the arrows of the quiver $\overline{A}_n$ on which $X$ is nonzero coincide with the arrow support of $\sigma(X)$.
\end{remark}

\begin{figure}
    \begin{tikzpicture}
        \draw (2,0) circle (2);
        \draw[fill=gray] (1,0) circle (0.25);
        \draw[fill=gray] (2,0) circle (0.25);
        \draw[fill=gray] (3,0) circle (0.25);

        \filldraw (2,2) circle (2.4pt);
        \filldraw (1.25,0) circle (2.4pt);
        \filldraw (2.25,0) circle (2.4pt);
        \filldraw (3.25,0) circle (2.4pt);
        \filldraw (4,0) circle (2.4pt);
        \filldraw (0,0) circle (2.4pt);
        \filldraw (2,-2) circle (2.4pt);
        
        \draw [smooth,dashed] (2,2) [out = -150,in=90] to (0.5,0) [out = -90, in = 150] to (2,-2);

        \draw [smooth,dashed] (2,2) [out = -120,in=90] to (1.5,0) [out = -90, in = 120] to (2,-2);

        \draw [smooth,dashed] (2,2) [out = -60,in=90] to (2.5,0) [out = -90, in = 60] to (2,-2);

        \draw [smooth,dashed] (2,2) [out = -30,in=90] to (3.5,0) [out = -90, in = 30] to (2,-2);
    \end{tikzpicture}
    \caption{A marked unpunctured surface which allows for the arc model of $\brick(\Pi(A_4,\{2,3\}))$ to be realized as a special case of Baur and Coelho Sim\~{o}es's surface model. See Remark~\ref{rem:baur_simoes}}\label{fig:BCS}
\end{figure}
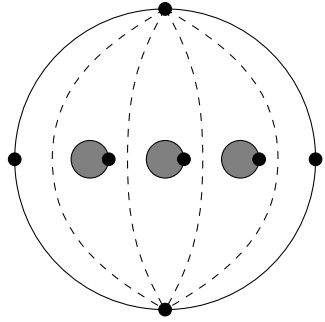

\begin{remark}\label{rem:baur_simoes}
    In \cite{BCS}, Baur and Coelho Sim\~{o}es describe the string modules over an arbitrary gentle algebra using ``homotopy classes of permissible arcs'' on partially-triangulated marked unpunctured surfaces (with boundary). We observe that the arc model for $\brick(\Pi(A_n,[2,n-1]_\mathbb{Z}))$ can be seen as a special case of Baur and Coelho Sim\~{o}es's model as follows. Start with a disk $D^2 \subseteq \mathbb{R}^2$ of radius $n/2$ centered at $(n/2,0)$. For each node $i \in [n]\setminus\{0,n\}$, we add a boundary component of radius 1/4 centered at the point $(i,0)$. Place marked points at the points $(n/2,n/2)$, $(n/2,-n/2)$, $(0,0)$, $(n,0)$, and $(i+1/4,0)$ for $i \in [n] \setminus \{0,n\}$. We then partially triangulate the resulting marked surface by drawing $n$ distinct edges between the marked points $(n/2,n/2)$ and $(n/2,-n/2)$, one passing through each point $(i+1/2,0)$ for $i \in [n] \setminus \{n\}$. See Figure~\ref{fig:BCS} for an example with $n = 4$. The ``tiling algebra'' of this triangulated marked surface, as defined in \cite[Definition~2.1]{BCS} is then precisely $\Pi(A_n,[2,n-1]_\mathbb{Z})$. Moreover, each ``homotopy class of permissible arcs'' can be seen as an arc connecting marked points along the $x$-axis in an analogous way to combinatorial equivalence classes of arcs.
\end{remark}

\section{Morphisms between bricks}\label{sec:hom_spaces}

In this section, we show how the arc model for $\brick(\Pi(A_n,S))$ can be used to describe the Hom-spaces between bricks. More precisely, recall from Section~\ref{sec:background} that $\Pi(A_n,[2,n-1])$ is a ``gentle algebra'' with no ``bands''. As such, an explicit basis for the Hom-space between two bricks is given in \cite[Section~2]{CB_string}. This basis is given an explicit interpretation in terms of arcs in \cite[Proposition~4.8]{BCZ}. It is also shown in \cite[Section~3.1]{mizuno2} that the basis is still valid over the preprojective algebra $\Pi(A_n)$, which can also be seen from the fact that $F_{[2,n-1]}^\emptyset$ is fully faithful. We recall the definition of this basis in Proposition~\ref{prop:quotients}(3). We then relate the cardinality of this basis to the set of intersection points between the corresponding arcs in Proposition~\ref{prop:hom_space}.

We begin with the following definitions.

\begin{definition}\label{def:quotient_arc}
    Let $\gamma_1,\gamma_2 \in \arc(n)$. For $t \in \{1,2\}$, denote $\wrd(\gamma_t) = \left(l(\gamma_t),s_{l(\gamma_t)+1}^t\cdots s_{r(\gamma_t)}^t\right)$.
    \begin{enumerate}
        \item We say that $\gamma_2$ is a \emph{restriction} of $\gamma_1$ if (i) $l(\gamma_1) \leq l(\gamma_2) < r(\gamma_2) \leq r(\gamma_1)$, and (ii) for every node $i$ in the arrow support of $\gamma_2$, one has $s_i^1 = s_i^2$. That is, the support of $\gamma_2$ is contained in that of $\gamma_1$ and the arcs pass on the same side of every node in their common arrow support.
        \item We say that $\gamma_2$ is a \emph{quotient arc} of $\gamma_1$ if it is a restriction of $\gamma_1$ which satisfies both of the following.
        \begin{enumerate}
            \item If $\gamma_1$ and $\gamma_2$ do not share a left endpoint, then $s_{l(\gamma_2)}^1 = \text{o}$; that is, $\gamma_1$ passes above the left endpoint of $\gamma_2$.
            \item If $\gamma_1$ and $\gamma_2$ do not share a right endpoint, then $s_{r(\gamma_2)}^1 = \text{u}$; that is, $\gamma_1$ passes below the right endpoint of $\gamma_2$.
        \end{enumerate}
        \item We say that $\gamma_2$ is a \emph{submodule arc} of $\gamma_1$ if it is a restriction of $\gamma_1$ which satisfies the following.
        \begin{enumerate}
            \item If $\gamma_1$ and $\gamma_2$ do not share a left endpoint, then $s_{l(\gamma_2)}^1 = \text{u}$; that is, $\gamma_1$ passes below the left endpoint of $\gamma_2$.
            \item If $\gamma_1$ and $\gamma_2$ do not share a right endpoint, then $s_{r(\gamma_2)}^1 = \text{o}$; that is, $\gamma_1$ passes above the right endpoint of $\gamma_2$.
        \end{enumerate}
    \end{enumerate}
\end{definition}

\begin{figure}
\begin{tikzpicture}[scale=0.9]
\draw [blue, very thick, <-,dotted]
(-0.3, -0.22+0.15) .. controls (-0.15, -0.5+0.15)and(0.15, -0.5+0.15) .. (0.3, -0.22+0.15);

\draw [blue, very thick, ->,dotted] 
(2-0.3, 0.22-0.15) .. controls (2-0.15, 0.5-0.15)and(2+0.15, 0.5-0.15) .. (2+0.3, 0.22-0.15);

\draw [blue, very thick, ->,dotted] 
(4-0.3, 0.22-0.15) .. controls (4-0.15, 0.5-0.15)and(4+0.15, 0.5-0.15) .. (4+0.3, 0.22-0.15);

\draw [blue, very thick,dotted] 
(-2,0) to [out = 45,in = 135]node[midway, above]{$\sigma(X)$}(1,0) to [out = -45,in = -135] (6,0);

;

\draw [orange, very thick]
(-4+4,0)..controls (-3.5+4, 0.3)and(-2.5+4,0.3)..(-2+4,0)
node[near end, above]{$\sigma(Y)$}
;

\draw [purple, very thick,dashed]
(-4+6,0)..controls (-3.5+6, -0.3)and(-2.5+6,-0.3)..(-2+6,0)
node[near end, below]{$\sigma(Z)$}
;
        
\filldraw  (0,0) circle (2.4pt)
 (-2,0) circle (2.4pt)
 (2,0) circle (2.4pt)
 (4,0) circle (2.4pt) 
 (6,0) circle (2.4pt)
 ;

\end{tikzpicture}
\caption{$\sigma(Y)$ (solid orange) is a quotient arc of $\sigma(X)$ (dotted blue). $\sigma(Z)$ (dashed purple) is a restriction of $\sigma(X)$, but is neither a quotient arc nor a submodule arc.}\label{fig:quotients}
\end{figure}
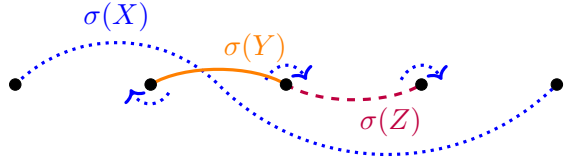

\begin{remark}
    Visually, we have that a restriction $\gamma_2$ of $\gamma_1$ is a quotient arc if and only if the small clockwise arrow of $\gamma_1$ at any non-shared endpoint of $\gamma_2$ points ``outward''. Likewise, we have that a restriction $\gamma_2$ of $\gamma_1$ is a submodule arc of $\gamma_1$ if and only if the small clockwise arrow of $\gamma_1$ at any non-shared endpoint of $\gamma_2$ points ``inward''. See Example~\ref{ex:quotients} and Figure~\ref{fig:quotients}.
\end{remark}

The following is contained in \cite[Section~3.1]{mizuno2} and \cite[Propositions~4.7 and~4.8]{BCZ}.

\begin{proposition}\label{prop:quotients} Let $S \subseteq [2,n-1]_\mathbb{Z}$.
    \begin{enumerate}
        \item Let $X \in \brick(\Pi(A_n,S))$ and let $Y \in \mods \Pi(A_n)$ be indecomposable. Then the following are equivalent.
        \begin{enumerate}
            \item $Y$ is a brick and $\sigma(Y)$ is a quotient arc of $\sigma(X)$.
            \item $Y$ is a quotient of $X$ in $\mods\Pi(A_n)$.
            \item $Y \in \mods\Pi(A_n,S)$ and is a quotient of $X$ in $\mods\Pi(A_n,S)$.
        \end{enumerate}
        Moreover, if these equivalent conditions hold then $$\dim\Hom_{\Pi(A_n,S)}(X,Y) = 1\qquad\text{and}\qquad\dim\Hom_{\Pi(A_n,S)}(Y,X) = 0.$$
        \item Let $X \in \brick(\Pi(A_n,S))$ and let $Y \in \mods \Pi(A_n)$ be indecomposable. Then the following are equivalent.
        \begin{enumerate}
            \item $Y$ is a brick and $\sigma(Y)$ is a submodule arc of $\sigma(X)$.
            \item $Y$ is a submodule of $X$ in $\mods \Pi(A_n)$.
            \item $Y \in \mods \Pi(A_n,S)$ and is a submodule of $X$ in $\mods \Pi(A_n,S)$.
        \end{enumerate}
        Moreover, if these equivalent conditions hold then $$\dim\Hom_{\Pi(A_n,S)}(X,Y) = 0\qquad\text{and}\qquad\dim\Hom_{\Pi(A_n,S)}(Y,X) = 1.$$
        \item Let $X, Y \in \brick(\Pi(A_n,S))$. Then the set of arcs $\gamma$ which are both quotient arcs of $\sigma(X)$ and submodule arcs of $\sigma(Y)$ form a basis of $\Hom_{\Pi(A_n,S)}(X,Y)$. The basis element corresponding to such an arc $\gamma$ is the composition of the quotient map $X \twoheadrightarrow \sigma^{-1}(\gamma)$ and the inclusion map $\sigma^{-1}(\gamma) \hookrightarrow Y$, each of which is unique up to scalar multiplication.
    \end{enumerate}
\end{proposition}

\begin{example}\label{ex:quotients}
Let $X = {\footnotesize\begin{matrix}2\\1 \ 3\\ \ \ \ 4\end{matrix}}, Y = 2 = S(2)$, and $Z = 3 = S(3)$. Then $\sigma(Y)$ is a quotient arc of $\sigma(X)$ and $\sigma(Z)$ is a restriction of $\sigma(X)$ which is neither a quotient arc nor a submodule arc, as shown in Figure~\ref{fig:quotients}. Moreover, we see that $\dim\Hom_{\Pi(A_n)}(X,Y) = 1$ and $\dim\Hom_{\Pi(A_n)}(X,Z) = 0$.
\end{example}

We now wish to relate the basis in Proposition~\ref{prop:quotients} to certain intersection points between the corresponding arcs. To do so, we need the following definitions.

\begin{definition}\label{def:crossing}
    Let $\gamma_1, \gamma_2 \in \arc(n)$, and choose representatives of their equivalences classes which intersect at a minimal number of points. Let $x \in \mathbb{R}$ be in the common closed support of $\gamma_1$ and $\gamma_2$, and suppose that $\gamma_1(x) = \gamma_2(x)$.
    \begin{enumerate}
        \item We say that $x$ is a \emph{shared (left or right) endpoint} of $\gamma_1$ and $\gamma_2$ if (left) $l(\gamma_1) = x = l(\gamma_2)$ or (right) $r(\gamma_1) = x = r(\gamma_2)$.
        \item We say that $x$ is a \emph{contested endpoint} of $\gamma_1$ and $\gamma_2$ if $r(\gamma_1) = x = l(\gamma_2)$ or $l(\gamma_1) = x = r(\gamma_2)$.
        \item We say that $x$ is a \emph{nontrivial crossing} of $\gamma_1$ and $\gamma_2$ if it is neither a shared nor a contested endpoint. Note that nontrivial crossings will necessarily be transverse by the minimality assumption.
        \item We say that $x$ is a \emph{Ext-crossing} of $\gamma_1$ and $\gamma_2$ if it is either a contested endpoint or a nontrivial crossing.
        \item We say that $x$ is a \emph{Hom-crossing} of $\gamma_1$ and $\gamma_2$ if it is either a shared endpoint or a nontrivial crossing. We say that a Hom-crossing $x$ is \emph{directed from $\gamma_1$ to $\gamma_2$} if one of the following holds.
        \begin{enumerate}
            \item $x$ is not a shared left endpoint and $\gamma_1$ is ``left locally above'' $\gamma_2$ near $x$. More precisely, for all $\varepsilon > 0$ there exists $0 < \delta < \varepsilon$ such that $\gamma_1(x-\delta) > \gamma_2(x-\delta)$.
            \item $x$ is a shared left endpoint and $\gamma_1$ is ``right locally below'' $\gamma_2$ near $x$. More precisely, for all $\varepsilon > 0$ there exists $0 < \delta < \varepsilon$ such that $\gamma_1(x+\delta) < \gamma_2(x + \delta)$.
        \end{enumerate}
    \end{enumerate}
\end{definition}

\begin{remark}
    Note that for $\gamma_1, \gamma_2 \in \arc(n)$, the number of crossings of each type only depends on the assumption that the chosen representatives intersect a minimal number of times, and not on the specific choice of representatives satisfying this assumption.
\end{remark}

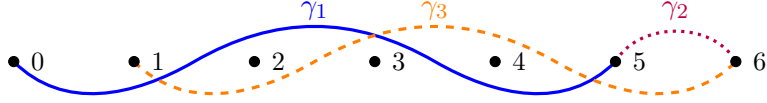
\begin{figure}
\begin{tikzpicture}[scale=0.8]

\draw[blue,very thick,smooth] (-4,0) [out = -45,in = -150] to (-1,0) [out = 30,in = 150] to (3,0) [out = -30,in = -135] to (6,0);

\draw[orange,dashed,very thick,smooth] (-2,0) [out = -45,in = -150] to (1,0) [out = 30,in = 150] to (5,0) [out = -30,in = -135] to (8,0);

\draw[purple,very thick,dotted] (6,0) [out = 60,in = 120] to (8,0);

\node at (1,0.85) {{\color{blue}$\gamma_1$}};
\node at (3,0.85) {{\color{orange}$\gamma_3$}};
\node at (7,0.85) {{\color{purple}$\gamma_2$}};

\node at (-3.6,0) {{\small 0}};
\node at (-1.6,0) {{\small 1}};
\node at (0.4,0) {{\small 2}};
\node at (2.4,0) {{\small 3}};
\node at (4.4,0) {{\small 4}};
\node at (6.4,0) {{\small 5}};
\node at (8.4,0) {{\small 6}};

\filldraw  (0,0) circle (2.4pt)
 (-2,0) circle (2.4pt)
 (-4,0) circle (2.4pt)
 (2,0) circle (2.4pt)
 (4,0) circle (2.4pt) 
 (6,0) circle (2.4pt)
 (8,0) circle (2.4pt)
 ;
\end{tikzpicture}
\caption{The arcs in Example~\ref{ex:crossings} $\gamma_1$ is drawn in solid blue, $\gamma_2$ in dotted purple, and $\gamma_3$ in dashed orange.}\label{fig:crossings}
\end{figure}

\begin{example}\label{ex:crossings}
    Let $\gamma_1$, $\gamma_2$, and $\gamma_3$ by the arcs in Figure~\ref{fig:crossings}. Then:
    \begin{enumerate}
        \item $\gamma_1$ and $\gamma_2$ have one Ext-crossing and zero Hom-crossings. The unique Ext-crossing is the contested endpoint 5.
        \item $\gamma_2$ and $\gamma_3$ have zero Ext-crossings and one Hom-crossing. The unique Hom-crossing is the shared (right) endpoint 6. This crossing is directed from $\gamma_2$ to $\gamma_3$.
        \item $\gamma_1$ and $\gamma_3$ have three nontrivial crossings, which count as both Hom-crossings and Ext-crossings. The crossings between the nodes 1 and 2 and between the nodes 4 and 5 are directed from $\gamma_3$ to $\gamma_1$, while the other crossing is directed from $\gamma_1$ to $\gamma_3$. Note also that the arcs can be perturbed so that the crossing directed from $\gamma_1$ to $\gamma_3$ occurs anywhere within the interval $(2,4) \subseteq \mathbb{R}$. On the other hand, the two crossings directed from $\gamma_3$ to $\gamma_1$ will always occur in the intervals $(1,2)$ and $(4,5)$, respectively.
    \end{enumerate}
\end{example}

Example~\ref{ex:crossings}(3) shows that the pair of nodes between which a crossing occurs is not invariant under combinatorial equivalence. We do, however, have the following.

\begin{lemma}\label{lem:crossing_exists}
    Let $\gamma_1, \gamma_2 \in \arc(n)$, and choose representatives of $\gamma_1$ and $\gamma_2$ which intersect a minimal number of times. For $t \in \{1,2\}$, denote $\wrd(\gamma_t) =: \left(l(\gamma_t),s^t_{l(\gamma_t)+1}\cdots s^t_{r(\gamma_t)}\right)$.
    \begin{enumerate}
        \item Suppose there is a Hom-crossing directed from $\gamma_1$ to $\gamma_2$ at some point $x$. Then there exist unique intergers $i<j$ in the common closed support of $\gamma_1$ and $\gamma_2$ which satisfy all of the following.
        \begin{enumerate}
            \item $i \leq x \leq j$.
            \item $s_k^1 = s_k^2$ for all $i < k < j$; that is, $\gamma_1$ and $\gamma_2$ pass on the same side of every node in the interval $(i,j)$.
            \item $s_{j}^1 \leq \textnormal{e} \leq s_{j}^2$; that is, either $j$ is a shared right endpoint or $\gamma_1$ must be below $\gamma_2$ at $j$.
            \item If $l(\gamma_1) \neq i$ then $s_i^1 = \textnormal{o}$; that is, either $i$ is the left endpoint of $\gamma_1$ or $\gamma_1$ passes above $i$.
            \item If $l(\gamma_2) \neq i$ then $s_i^2 = \textnormal{u}$; that is, either $i$ is the left endpoint of $\gamma_2$ or $\gamma_2$ passes below $i$.
        \end{enumerate}
        Moreover, $x$ is a nontrivial crossing if and only if $x \notin \{i,j\}$. Otherwise, $x$ is a shared endpoint.
    \item Suppose there exist integers $i < j$ in the common closed support of $\gamma_1$ and $\gamma_2$ which satisfy conditions (b)-(e) from (1). Then either $\gamma_1$ and $\gamma_2$ are combinatorially equivalent or there is a unique crossing between $\gamma_1$ and $\gamma_2$ in the interval $[i,j] \subseteq \mathbb{R}$. Moreover, in the second case, this crossing is a Hom-crossing directed from $\gamma_1$ to $\gamma_2$.
    \end{enumerate}
\end{lemma}

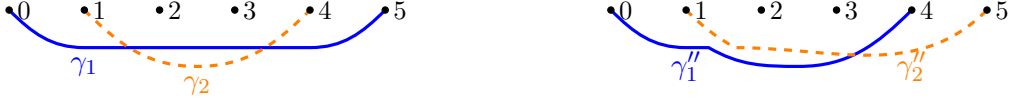
\begin{figure}
\begin{tikzpicture}[scale=0.5]

\draw [blue, very thick,smooth] 
(-2,0) [out = -45, in = 180]  to (0,-1) [out = 0,in = 180] to (6,-1)[out = 0, in = -135] to (8,0);

\node[color = blue] at (0,-1.5){$\gamma_1$};

\draw [orange, very thick,smooth,dashed] 
(0,0) [out = -45, in = 180]  to (3,-1.5) [out = 0,in = -135] to (6,0);

\node[color = orange] at (3,-2){$\gamma_2$};
        
\filldraw  (0,0) circle (2.4pt)
 (-2,0) circle (2.4pt)
 (2,0) circle (2.4pt)
 (4,0) circle (2.4pt) 
 (6,0) circle (2.4pt)
 (8,0) circle (2.4pt)
 ;

\node at (-1.6,0) {{\small 0}};
\node at (0.4,0) {{\small 1}};
\node at (2.4,0) {{\small 2}};
\node at (4.4,0) {{\small 3}};
\node at (6.4,0) {{\small 4}};
\node at (8.4,0) {{\small 5}};

\begin{scope}[shift = {(16,0)}]
\draw [blue, very thick,smooth] 
(-2,0) [out = -45, in = 180]  to (0,-1) [out = 0,in = 180] to (0.6,-1)[out = -30, in = 180] to (3,-1.5) [out = 0, in = -135] to (6,0);

\node[color = blue] at (0,-1.5){$\gamma_1''$};

\draw [orange, very thick,smooth,dashed] 
(0,0) [out = -45, in = 150]  to (1.3,-1) [out = 0,in = -135] to (8,0);

\node[color = orange] at (6,-1.5){$\gamma_2''$};
        
\filldraw  (0,0) circle (2.4pt)
 (-2,0) circle (2.4pt)
 (2,0) circle (2.4pt)
 (4,0) circle (2.4pt) 
 (6,0) circle (2.4pt)
 (8,0) circle (2.4pt)
 ;

\node at (-1.6,0) {{\small 0}};
\node at (0.4,0) {{\small 1}};
\node at (2.4,0) {{\small 2}};
\node at (4.4,0) {{\small 3}};
\node at (6.4,0) {{\small 4}};
\node at (8.4,0) {{\small 5}};
\end{scope}

\end{tikzpicture}
\caption{An illustration of the proof of Lemma~\ref{lem:crossing_exists}. The fact that $\gamma_1(4) < 0$ implies that the two crossings between $\gamma_1$ and $\gamma_2$ are unnecessary.}\label{fig:resolve}
\end{figure}

\begin{proof}
    (1) Suppose first that $x$ is a nontrivial crossing directed from $\gamma_1$ to $\gamma_2$. By perturbing $\gamma_1$ and $\gamma_2$ if necessary, we can assume that $x \notin \mathbb{Z}$. By definition, there exists $\delta > 0$ such that $\gamma_1(x-\delta) > \gamma_2(x-\delta)$ and the interval $(x-\delta,x+\delta)$ does not contain an integer. Moreover, the minimality assumption on $\gamma_1$ and $\gamma_2$ implies that the intersection at $x$ is transverse. In particular, we can assume that $x$ is the only crossing between $\gamma_1$ and $\gamma_2$ in the interval $(\lfloor x-\delta \rfloor, \lceil x+\delta \rceil)$, where $\lfloor -\rfloor$ and $\lceil - \rceil$ denote the floor and ceiling functions, respectively.
    
    We now define two new arcs $\gamma_1'$ and $\gamma_2'$ such that
    $$\gamma_1'(t) = \begin{cases} \gamma_1(t) & t \leq x\\ \gamma_2(t) & t \geq x,\end{cases}\qquad\qquad\gamma_2'(t) = \begin{cases} \gamma_2(t) & t \leq x\\ \gamma_1(t) & t \geq x.\end{cases}$$
    The combinatorial data determined by these arcs are    
        \begin{eqnarray*}
            \wrd(\gamma_1') &=& \left(l(\gamma_1),s^1_{l(\gamma_1)+1}\cdots s^1_{\lfloor x\rfloor}s^2_{\lceil x \rceil} \cdots s^2_{r(\gamma_2)}\right),\\
            \wrd(\gamma_2') &=& \left(l(\gamma_2),s^2_{l(\gamma_2)+1}\cdots s^2_{\lfloor x\rfloor}s^1_{\lceil x \rceil} \cdots s^1_{r(\gamma_2)}\right).
        \end{eqnarray*}
    Note that the crossing between $\gamma_1'$ and $\gamma_2'$ at $x$ is not transverse. Thus we can perturb $\gamma_1'$ and $\gamma_2'$ into arcs $\gamma_1''$ and $\gamma_2''$ which agree with $\gamma_1'$ and $\gamma_2'$ outside the interval $(\lfloor x-\delta \rfloor, \lceil x+\delta \rceil)$ and do not intersect within this interval. In particular, $\gamma_1''$ is combinatorially equivalent to $\gamma_1'$ and likewise for $\gamma_2''$. Moreover, $\gamma_1''$ and $\gamma_2''$ have fewer nontrivial crossings than $\gamma_1$ and $\gamma_2$.

    Now suppose for a contradiction that there do not exist $i < j$ satisfying the conditions (a)-(e). We then have two cases to consider.

    Suppose first that condition (c) must fail. Then for all nodes $j$ with $x < j \leq \min\{r(\gamma_1),r(\gamma_2)\}$, we must have $s_j^1 \not< s_j^2$. If $s_j^1 = s_j^2$ for all nodes $j \in (x,\min\{r(\gamma_1),r(\gamma_2)\}]$, then $\gamma_1$ is combinatorially equivalent to $\gamma_1''$, and likewise for $\gamma_2$. This contradicts the fact that the chosen representatives of $\gamma_1$ and $\gamma_2$ intersect a minimal number of times. Otherwise, let $j \in (x, \min\{r(\gamma_1),r(\gamma_2)\}]$ be minimal such that $s_j^1 \neq s_j^2$. Then $s_j^1 > s_j^2$, and so $\gamma_2''(j) > \gamma_1''(j)$ by the assumption that (c) fails. On the other hand, we have $\gamma_2''(x) < \gamma_1''(x)$ by construction. But this means there is an intersection between $\gamma_1''$ and $\gamma_2''$, and hence also between $\gamma_1$ and $\gamma_2$, somewhere in the interval $(x,j)$.The minimality of $j$ then implies that $\gamma_1''$ and $\gamma_1$ pass on the same side of every node in $(x,j)$, and likewise for $\gamma_2''$ and $\gamma_2$. Thus the crossings between $\gamma_1$ and $\gamma_2$ at $x$ and in $(x,j)$ can be resolved, a contradiction. See Figure~\ref{fig:resolve} for an example.

    The other case is that at least one of conditions (d) and (e) must fail. This leads to a similar contradiction as in the previous case.

    This concludes the proof in the case that $x$ is a nontrivial crossing. Suppose next that $x$ is a shared left endpoint and take $i = x = l(\gamma_1) = l(\gamma_2)$. Then $\gamma_1$ is right locally below $\gamma_2$ near $x$, and so $s_{i + 1}^1 \leq s_{i+1}^2$. It follows that either $s_{j}^1 = s_j^2$ for all $i < j \leq \min\{r(\gamma_1),r(\gamma_2)\}$ or there exists some minimal $i < j \leq \min\{r(\gamma_1),r(\gamma_2)\}$ such that $s_{j}^1 < s_{j}^2$. The proof when $x$ is a shared right endpoint is similar.

    (2) Suppose conditions (b)-(e) are satisfied for some $i < j$. If $i$ and $j$ are both shared endpoints between $\gamma_1$ and $\gamma_2$, then $\wrd(\gamma_1) = \wrd(\gamma_2)$ and we are done. Thus assume this is not the case. Then $\gamma_1(i) \geq \gamma_2(i)$ and $\gamma_1(j) \leq \gamma_2(j)$, with at least one of these inequalities being strict. This implies that $\gamma_1$ and $\gamma_2$ must cross somewhere in the interval $[i,j]$. To show that this crossing is unique, suppose there exist $i \leq x < y \leq j$ such that $\gamma_1$ and $\gamma_2$ cross at both $x$ and $y$. We suppose that $y \neq j$, the case where $i \neq x$ being similar. Furthermore, we can assume without loss of generality that $\gamma_1$ and $\gamma_2$ do not cross in the interval $(x,y)$. Then, by assumption, $\gamma_1$ and $\gamma_2$ pass on the same side of every node in $(x,y)$. Using a similar argument as in the proof of (1), this contradicts that the chosen representatives of $\gamma_1$ and $\gamma_2$ intersect a minimal number of times.
\end{proof}

We now combine Proposition~\ref{prop:quotients} and Lemma~\ref{lem:crossing_exists} to deduce the main result of this section.

\begin{proposition}\label{prop:hom_space}
    Let $S \subseteq [2,n-1]_\mathbb{Z}$ and let $X, Y \in \brick(\Pi(A_n,S))$. Choose representatives of the corresponding arcs $\sigma(X)$ and $\sigma(Y)$ which intersect a minimal number of times. Then $\dim\Hom_{\Pi(A_n,S)}(X,Y)$ is equal to the number of Hom-crossings directed from $\sigma(X)$ to $\sigma(Y)$.
\end{proposition}

\begin{proof}
    By Proposition~\ref{prop:quotients}, it suffices to show that the number of Hom-crossings directed from $\sigma(X)$ to $\sigma(Y)$ coincides with the number of (combinatorial equivalence classes of) arcs which are both quotient arcs of $\sigma(X)$ and submodule arcs of $\sigma(Y)$.

    Suppose first that $X = Y$. Then $\sigma(X) = \sigma(Y)$ is the only arc which is both a quotient arc of $\sigma(X)$ and a submodule arc of $\sigma(Y)$. Moreover, $\sigma(X)$ and $\sigma(Y)$ have two shared endpoints, one directed from $\sigma(X)$ to $\sigma(Y)$ as a Hom-crossing and one directed from $\sigma(Y)$ to $\sigma(X)$ as a Hom-crossing. This proves the result in this case.
    
    Now suppose $X \neq Y$. First let $\gamma \in \arc(n)$ be a quotient arc of $\sigma(X)$ and a submodule arc of $\sigma(Y)$. Then $i = l(\gamma)$ and $j = r(\gamma)$ must satisfy conditions (b)-(e) of Lemma~\ref{lem:crossing_exists}(1) for $\gamma_1 = \sigma(X)$ and $\gamma_2 = \sigma(Y)$. Thus, by Lemma~\ref{lem:crossing_exists}(2), there exists a crossing between $\gamma_1$ and $\gamma_2$ in the interval $[i,j]$, and this crossing is directed from $\sigma(X)$ to $\sigma(Y)$ as a Hom-crossing.

    Conversely, let $x$ be a Hom-crossing directed from $\sigma(X)$ to $\sigma(Y)$. Then there exist unique intergers $i < j$ which satisfy conditions (a)-(e) in Lemma~\ref{lem:crossing_exists}(1) for $\gamma_1 = \sigma(X)$ and $\gamma_2 = \sigma(Y)$. Now let $\gamma \in \arc(n)$ be an arc with endpoints $l(\gamma) = i$ and $r(\gamma) = j$ which passes on the same side of every node in $(i,j)$ as the arcs $\sigma(X)$ and $\sigma(Y)$. Then $\gamma$ is both a quotient arc of $\sigma(X)$ and a submodule arc of $\sigma(Y)$.

    Now note that the processes in the previous two paragraphs are inverse to one another. This proves the result.
\end{proof}

In particular, we note that Proposition~\ref{prop:hom_space} includes the characterization of semibricks over $\Pi(A_n,S)$ in terms of \emph{noncrossing arc diagrams} given in \cite[Section~4]{BCZ}. (Only the algebra $\Pi(A_n,[2,n-1])$ is explicitly considered in \cite{BCZ}. See also \cite[Theorem~3.13]{mizuno2} for an explicit proof of this characterization over $\Pi(A_n)$.) We recall the definition of a noncrossing arc diagram from \cite{reading}.

\begin{definition}\label{def:noncrossing}
    Let $\mathcal{A} \subseteq \arc(n)$ be a set of arcs. We say that $\mathcal{A}$ is a \emph{noncrossing arc diagram} if no distinct $\gamma \neq \rho \in \mathcal{A}$ have a hom-crossing. More precisely, there exist representatives of the combinatorial equivalence classes of $\gamma$ and $\rho$ which either do not intersect or intersect only at a contested endpoint. We denote by $\arc_{nc}(n)$ the set of noncrossing arc diagrams (up to combinatorial equivalence) on $n$ nodes.
\end{definition}

\begin{remark}
    While it does not follow from the definition \emph{a priori}, Reading shows in \cite[Proposition~3.2]{reading} that $\mathcal{A}$ is a noncrossing arc diagram if and only if there are \emph{fixed} representatives of the arcs in $\mathcal{A}$ which only have intersections at contested endpoints. This justifies the use of the name ``noncrossing arc diagram'' instead of ``pairwise noncrossing arc diagram''.
\end{remark}

The semibricks over $\Pi(A_n,S)$ are then characterized as follows.

\begin{corollary}\label{cor:noncrossing}\cite[Section~4]{BCZ}\cite[Theorem~3.13]{mizuno2}
    Let $S \subseteq [2,n-1]_\mathbb{Z}$. Then the map $\sigma$ induces a bijection $\sbrick(\Pi(A_n,S)) \rightarrow \arc_{nc}(n)$.
\end{corollary}

Another consequence of Proposition~\ref{prop:hom_space} is the following.

\begin{corollary}\label{cor:hom}
    Let $S \subseteq [2,n-1]_\mathbb{Z}$ and let $X, Y \in \brick(\Pi(A_n,S))$. Then $$\left|\dim\Hom_{\Pi(A_n,S)}(X,Y) - \dim\Hom_{\Pi(A_n,S)}(Y,X)\right| \leq 1.$$
\end{corollary}

\begin{proof}
    Choose representatives of the arcs $\sigma(X)$ and $\sigma(Y)$ which intersect a minimal number of times. Let $x_1 < \cdots < x_k$ be the set of Hom-crossings between these representatives. Then for $1 \leq i < k$, the minimality of the chosen representatives implies that if $x_i$ is directed from $\sigma(X)$ to $\sigma(Y)$, then $x_{i+1}$ is directed from $\sigma(Y)$ to $\sigma(X)$, and vice versa. The result thus follows from Proposition~\ref{prop:hom_space}.
\end{proof}

\section{Extensions between bricks}\label{sec:ext}

We now turn our attention to using the arc model to compute a basis for the first Ext-space between bricks. Unlike the Hom-spaces, the Ext-spaces will generally depend on the choice of $S$ when working over $\Pi(A_n,S)$. See Example~\ref{ex:extension3} for an example. However, we recall from Section~\ref{sec:preproj} that given $S \subseteq T \subseteq [2,n-1]_\mathbb{Z}$ and $M, N \in \mods \Pi(A_n,T)$, the functor $F_T^S$ allows us to consider $\Ext^1_{\Pi(A_n,T)}(M,N)$ as a subspace of $\Ext^1_{\Pi(A_n,S)}(M,N)$.

To start, we recall the following formula of Crawley-Boevey.

\begin{proposition}\cite{CB}\label{prop:CB}
    Define a bilinear form $(-,-):\mathbb{R}^n \rightarrow \mathbb{R}$ by
    $$
        (\alpha,\beta) = \sum_{i = 1}^n 2\alpha_i\beta_i - \sum_{i = 1}^{n-1} (\alpha_i\beta_{i + 1} + \alpha_{i+1}\beta_i)
    $$
    Then for all $M, N \in \mods \Pi(A_n)$, one has
    $$(\undim M,\undim N) = \dim\Hom_{\Pi(A_n)}(M,N) + \dim\Hom_{\Pi(A_n)}(N,M) - \dim\Ext^1_{\Pi(A_n)}(M,N).$$
\end{proposition}

We note that Crawley-Boevey proves a version of Proposition~\ref{prop:CB} for the preprojective algebra of any quiver, but we have specified the result to $\Pi(A_n)$ for simplicity. Moreover, when $M$ and $N$ are bricks we can simplify Proposition~\ref{prop:CB} and use it to compute the dimension of $\Ext^1$ as follows.

\begin{corollary}\label{cor:CB}
    Let $X, Y \in \brick(\Pi(A_n))$. Choose representatives of the arcs $\sigma(X)$ and $\sigma(Y)$ which intersect a minimal number of times. Denote by $\mathrm{ce}(X,Y)$ the number of contested endpoints between $\sigma(X)$ and $\sigma(Y)$, by $\mathrm{se}(X,Y)$ the number of shared endpoints between $\sigma(X)$ and $\sigma(Y)$, and by $\mathrm{nc}(X,Y)$ the number of nontrivial crossings between $X$ and $Y$. Then the following hold.
    \begin{enumerate}
        \item $(\undim X,\undim Y) = \mathrm{se}(X,Y) - \mathrm{ce}(X,Y)$.
        \item $\dim\Ext^1_{\Pi(A_n)}(X,Y) = \mathrm{ce}(X,Y) + \mathrm{nc}(X,Y)$, which is precisely the number of Ext-crossings between $\sigma(X)$ and $\sigma(Y)$.
    \end{enumerate}
\end{corollary}

\begin{proof}
    (1) For readability, we append 0 to each end of $\alpha:= \undim X$ and $\beta:= \undim Y$; that is, we denote $\alpha_0 = 0 = \alpha_{n+1}$ and $\alpha_i = \dim X(i)$ for all other $i$ and likewise for $\beta$. We also denote $l_X:= l(\sigma(X))$, and likewise for $r_X, l_Y$, and  $r_Y$. Since $(-,-)$ is symmetric, we assume without loss of generality that $l_X \leq l_Y$.

    Suppose first that $r_X \leq l_Y$. This means the supports of $X$ and $Y$ have an empty intersection. Then $\sigma(X)$ and $\sigma(Y)$ have no shared endpoints and there are two possibilities. If $r_X = l_Y$ is a contested endpoint, then $\dim X(r_X-1) = 1 = \dim Y(r_X)$ and $(\undim X,\undim Y) = -1$. Otherwise, $X$ and $Y$ have no contested endpoints and $(\undim X,\undim Y) = 0$. This concludes the proof in this case.
    
    Suppose from now on that $l_Y < r_X$, and so in particular $\sigma(X)$ and $\sigma(Y)$ have no contested endpoints. We next consider the case where $Y$ is a simple module; i.e., where $l_Y = r_Y - 1$. In this case, we have $(\undim X,\undim Y) = 2\alpha_{r_Y} - \alpha_{r_Y+1} - \alpha_{r_Y-1}$, where $\alpha_{r_Y} \neq 0$ by assumption. Now the condition $\alpha_{r_Y+1} = 0$ is equivalent to $r_X = r_Y$ and the condition $\alpha_{r_Y-1} = 0$ is equivalent to $l_X = l_Y$. Thus the result holds in this case. The case where $X$ is a simple module is similar.
    
    Suppose from now on that neither $X$ nor $Y$ are simple. Again since $(-,-)$ is symmetric, we then have two cases left to consider: (i) $l_X \leq l_Y < r_Y \leq r_X$ and (ii) $l_X < l_Y < r_X < r_Y$.

    Suppose we are in case (i), and denote $z_i := (2\alpha_i-\alpha_{i+1}-\alpha_{i-1})\beta_i$. Now for $i \in [n]\setminus\{0\}$, the assumption of case (i) implies that if $\alpha_i = 0$ then $\beta_i = 0$. Moreover, since $X$ is not simple, we have that if $\alpha_i \neq 0$ then at least one of $\alpha_{i-1}$ and $\alpha_{i+1}$ is also nonzero. We conclude that $z_i \in \{0,1\}$ for all $i$. As $(\undim X,\undim Y) = \sum_{i = 1}^n z_i$, it remains only to show that the number of times $z_i$ takes the value 1 coincides with $\mathrm{se}(X,Y)$. Indeed, suppose $z_i = 1$, so in particular $\beta_i \neq 0$. If $\alpha_{i+1} = 0$, then we have $i = r_X = r_Y$, and this is a shared endpoint. Likewise, if $\alpha_{i-1} = 0$, then we have $i-1 = l_X = l_Y$ and this is a shared endpoint. This concludes the proof in this case.

    Finally, we suppose we are in case (ii). For $1 \leq i \leq n$, denote $w_i = 2\alpha_i\beta_i-\alpha_i\beta_{i+1}-\alpha_{i+1}\beta_i$. We then compute directly that $\alpha_i\beta_i = 1$ for $l_Y < i \leq r_X$ and is otherwise equal to 0, that $\alpha_{i+1}\beta_{i} = 1$ for $l_Y + 1 < i \leq r_X$ and is otherwise equal to 0, and that $\alpha_i\beta_{i+1} = 1$ for $l_Y < i \leq r_X + 1$ and is otherwise equal to 0. Thus we have $$w_i = \begin{cases} 1 & i = l_Y+1\\-1 & i = r_X + 1\\0 & \text{otherwise}.\end{cases}$$
    We conclude that $(\undim X,\undim Y) = 0$, as desired.

    (2) This follows immediately from combining (1) with Propositions~\ref{prop:hom_space} and~\ref{prop:CB}.
\end{proof}

Our next task is to construct an explicit basis of $\Ext^1_{\Pi(A_n)}(X,Y)$ for $X, Y \in \brick(\Pi(A_n))$. In doing so, we will also give a description for $\Ext^1_{\Pi(A_n,S)}(X,Y)$ for $S \subseteq [2,n-1]_\mathbb{Z}$. In particular, this will build off of the basis for $\Ext^1_{\Pi(A_n,[2,n-1]_\mathbb{Z})}(X,Y)$ described in \cite[Section~6.2]{BaH} (where it was deduced from \cite[Theorem~6.5]{BDMTY} and \cite[Theorem~A]{CPS}).

We first establish some notations via a series of examples. We will identify the following quiver $Q$ and representation $M_Q$ over $\overline{A_5}$:
$$Q = 
    \begin{tikzcd}
        & 2 \arrow[r]\arrow[dl]\arrow[dr] & 3 & 4\arrow[l]\\
        1 \arrow[r] & 2 & 3\arrow[l]\arrow[r] &4
    \end{tikzcd}
    \qquad\qquad
    M_Q = \begin{tikzcd} K\arrow[r,shift left,"\text{\tiny $\begin{bmatrix}0\\1\end{bmatrix}$}"] & K^2\arrow[l,shift left,"\text{\tiny $\begin{bmatrix}1 \ 0\end{bmatrix}$}"]\arrow[r,shift left,"\text{\tiny $\begin{bmatrix}1 \ 0\\1 \ 0\end{bmatrix}$}"] & K^2\arrow[l,shift left,"\text{\tiny $\begin{bmatrix}0 \ 0\\0 \ 1\end{bmatrix}$}"]\arrow[r,shift left,"\text{\tiny $\begin{bmatrix}0 \ 0\\ 0 \ 1\end{bmatrix}$}"] & K^2\arrow[l,shift left,"\text{\tiny $\begin{bmatrix}1 \ 0\\ 0\ 0\end{bmatrix}$}"]\arrow[r,shift left] & 0.\arrow[l,shift left]\end{tikzcd}
$$
Note that, in this case, we have that $M_Q \in \mods \Pi(A_5,\{3,4\})$ and that $M_Q \notin \Pi(A_5,\{2\})$. If in addition we wish to indicate that there is a (possibly empty) quiver of type $A_k$ ``attached'' to one end of our diagram, but that the directions of the arrows within this diagram are not relevant, we may write
$$Q = 
    \begin{tikzcd}
        & 2 \arrow[r,dashed,no head]\arrow[dl]\arrow[dr] & Q'\\
        1 \arrow[r] & 2 & 3\arrow[l]\arrow[r,dashed] &Q''.
    \end{tikzcd}
$$
Here, the quiver $Q'$ may be empty, a vertex 3, a pair of vertices 3 and 4 with an arrow either $3 \rightarrow 4$ or $4 \rightarrow 3$, etc., and similarly for $Q''$. The dashed arrow oriented from 3 to $Q''$ indicates that if $Q''$ is nonempty, then there must be an arrow $3\rightarrow 4$. The dashed arrow between 2 and $Q'$ indicates that if $Q'$ is nonempty, then there may either be an arrow $2\rightarrow 3$ or $2 \leftarrow 3$. In such a diagram, we assume that the only arrow connecting a vertex in $Q'$ to one outside of $Q'$ lies along the dotted line connecting 2 and $Q'$, and similarly for $Q''$. In particular, there are no arrows between $Q'$ and $Q''$.

Similarly, for $X \in \brick(\Pi(A_n))$, recall from Proposition~\ref{prop:arcBricks} that we can visualize $X$ as a quiver of type $A_n$ by drawing a vertex labeled $i$ if $X(i) \neq 0$ and drawing an arrow $X(i) \rightarrow X(i+1)$ (resp.\ $X(i) \leftarrow X(i+1))$ if $X(a_i) \neq 0$ (resp.\ if $X(a_i^*) \neq 0$). For example, suppose $\sigma(X) = \gamma_1$ as shown in Figure~\ref{fig:crossings}. Then we associate $X$ with the quiver
$$\begin{tikzcd} 1 \arrow[r] & 2 & 3\arrow[l] & 4\arrow[l]\arrow[r] & 5. \end{tikzcd}$$ If we are concerned primarily with the behavior of this brick between a certain subset of vertices, say $\{3,4,5\}$, we may abbreviate the quiver as
$$\begin{tikzcd}\textnormal{left}(X) \arrow[r,dashed,no head] & 3 & 4\arrow[l]\arrow[r] & 5 \arrow[r,dashed,no head] & \textnormal{right}(X),\end{tikzcd}$$
where here $\textnormal{right}(X)$ is empty and $\textnormal{left}(X) = 1\rightarrow 2 \leftarrow$.

We now associate short exact sequences to the various types of Ext-crossings between arcs. Let $S \subseteq [2,n-1]_\mathbb{Z}$, and let $X, Y \in \brick(\Pi(A_n,S))$.

Suppose first that $\sigma(X)$ and $\sigma(Y)$ have a contested endpoint $x$. If $l(\sigma(X)) = x = r(\sigma(Y))$, denote
    $$\begin{tikzcd}[column sep = small] Q = &\textnormal{left}(Y) \arrow[r,dashed,no head] & r(\sigma(Y)) & l(\sigma(X)) + 1 \arrow[l]\arrow[r,dashed,no head]& \textnormal{right}(X).\end{tikzcd}$$
    Otherwise, denote
    $$\begin{tikzcd}[column sep = small] Q = &\textnormal{left}(X) \arrow[r,dashed,no head] & r(\sigma(X))\arrow[r] & l(\sigma(Y)) + 1 \arrow[r,dashed,no head]& \textnormal{right}(Y).\end{tikzcd}$$
    Now let
\begin{equation}\label{extension1}\eta_1(x) = \qquad 0 \rightarrow Y \rightarrow M_Q \rightarrow X \rightarrow 0,\end{equation}
    with the maps at each vertex given by the identity whenever both vector spaces are nonzero. 

\begin{proposition}\label{prop:extension1}
    Let $\eta_1(x)$ be as in Equation~\eqref{extension1}. Then $\eta_1(x)$ is a nonsplit short exact sequence in $\Ext^1_{\Pi(A_n,S)}(X,Y)$. Moreover, recalling that we are in the case where $\sigma(X)$ and $\sigma(Y)$ have a contested endpoint, this short exact sequence spans $\Ext^1_{\Pi(A_n,S)}(X,Y)$.
\end{proposition}

\begin{proof}
    The result is shown for $S = [2,n-1]_\mathbb{Z}$ in \cite[Theorem~8.5]{BDMTY}. (One can also see that the sequence is well-defined and does not split for any $S$ since $M_Q$ is a brick.) Thus for an arbitrary $S \subseteq [2,n-1]_\mathbb{Z}$, we have
    $$0 \neq \eta_1(x) \in \Ext^1_{\Pi(A_n,[2,n-1]_\mathbb{Z})}(X,Y) \subseteq \Ext^1_{\Pi(A_n,S)}(X,Y) \subseteq \Ext^1_{\Pi(A_n)}(X,Y).$$
    Moreover, if $\sigma(X)$ and $\sigma(Y)$ have a contested endpoint, then they cannot have a nontrivial crossing, and so $\dim\Ext^1_{\Pi(A_n)}(X,Y) = 1$ by Corollary~\ref{cor:CB}. This proves the result.
\end{proof}

\begin{figure}
\begin{tikzpicture}[scale=0.8]

\draw[blue,very thick,smooth] (-4,0) [out = -45,in = 180] to (-2.5,-0.75) [out = 0,in = 180] to (0.5,0.75) [out = 0,in = 135] to (2,0);

\draw[blue,very thick,smooth] (2,0) [out = 45,in = 180] to (4,0.75) [out = 0,in = 145] to (6,0);

\draw[orange,dashed,very thick,smooth] (-4,0) [out = -30,in = 180] to (-2.5,-0.5) [out = 0,in = 180] to (0.5,1) [out = 0,in = 180] to (2,1) [out = 0, in = 180] to (4,1) [out = 0,in = 120] to (6,0);

\draw[purple,dotted,very thick,smooth] (-4,0) [out = -60,in = 180] to (-2.5,-1) [out = 0,in = 180] to (0.5,0.5) [out = 0,in = 180] to (2,-0.75) [out = 0, in = 180] to (3.75,0.5) [out = 0,in = 180] to (6,0);

\node at (0,1.2) {{\color{orange}$\gamma_U$}};
\node at (-1,-0.85) {{\color{purple}$\gamma_D$}};

\node at (-4.4,0) {{\small 0}};
\node at (-2.4,0) {{\small 1}};
\node at (0.4,0) {{\small 2}};
\node at (2.4,0) {{\small 3}};
\node at (4.4,0) {{\small 4}};
\node at (6.4,0) {{\small 5}};

\filldraw  (0,0) circle (2.4pt)
 (-2,0) circle (2.4pt)
 (-4,0) circle (2.4pt)
 (2,0) circle (2.4pt)
 (4,0) circle (2.4pt) 
 (6,0) circle (2.4pt)
 ;

 \draw [purple, very thick, ->,dotted] 
(2-0.3, 0.22+0.2) .. controls (2-0.15, 0.5+0.2)and(2+0.15, 0.5+0.2) .. (2+0.3, 0.22+0.2);

 \draw [orange, very thick, ->,dashed] 
(2+0.3, -0.22) .. controls (2+0.15, -0.5)and(2-0.15, -0.5) .. (2-0.3, -0.22);
\end{tikzpicture}
\caption{The arcs in Example~\ref{ex:extension1}. $\gamma_L$ is the solid blue arc with endpoints 0 and 3, and $\gamma_R$ is the solid blue arc with endpoints 3 and 5.}\label{fig:extension1}
\end{figure}
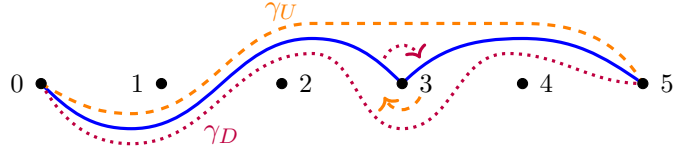

\begin{example}\label{ex:extension1}
    Let $\gamma_L$ and $\gamma_R$ be the solid blue arcs in Figure~\ref{fig:extension1} which have endpoints 0 and 3 and 3 and 5, respectively. Set $X = \sigma^{-1}(\gamma_L) = \text{\tiny $\begin{matrix}1 \ 3\\ 2\end{matrix}$}$ and $Y = \sigma^{-1}(\gamma_R) = \text{\tiny $\begin{matrix}5\\4\end{matrix}$}$. Then the term $M_Q$ in Equation~\eqref{extension1} is given by $M_Q = \text{\tiny $\begin{matrix}\ \ \ \ 5\\\ \ \ 4\\1 \ 3\\ 2\end{matrix}$}$.
    The arc $\sigma(M_Q)$ is precisely the dashed orange arc $\gamma_U$ in Figure~\ref{fig:extension1}, formed by gluing the arcs $\gamma_L$ and $\gamma_R$ at their contested endpoint 3 and perturbing upwards. In particular, the small clockwise arrow of $\gamma_U$ points away from $\gamma_R$ towards $\gamma_L$. If we instead perturbed the glued arc downwards to the dotted purple arc $\gamma_D$, the small clockwise arrow would point in the opposite direction and the corresponding short exact sequence would be an element of $\Ext^1_{\Pi(A_5,S)}(Y,X)$.
\end{example}

Suppose now that there is a nontrivial crossing $x$ directed from $\sigma(Y)$ to $\sigma(X)$. Let $i < j$ be the unique integers satisfying the properties (a)-(e) in Lemma~\ref{lem:crossing_exists} for $\gamma_1 = \sigma(Y)$ and $\gamma_2 = \sigma(X)$. Denote
    $$
    \begin{tikzcd}[column sep = small]
        Q= & \textnormal{left}(Y)  & i+1\arrow[l,dashed] \arrow[r,dashed,no head] & \cdots \arrow[r,dashed,no head] & j& \textnormal{right}(X)\arrow[l,dashed]\\
        Q'= & \textnormal{left}(X) \arrow[dashed,r] & i+1 \arrow[r,dashed,no head] & \cdots \arrow[r,dashed,no head] & j\arrow[r,dashed] & \textnormal{right}(Y),\\
    \end{tikzcd}
    $$
    where, for $i < k < j$, both $Q$ and $Q'$ contain an arrow $k \rightarrow k+1$ (resp.\ $k \leftarrow k+1$) if $\sigma(X)$ and $\sigma(Y)$ pass below (resp.\ above) the node $k$.
    Now let
    \begin{equation}\label{extension2}
        \eta_2(x) = \qquad 0 \rightarrow Y \rightarrow M_Q \oplus M_{Q'} \rightarrow X \rightarrow 0,
    \end{equation}
    The component map $Y \rightarrow M_Q$ is the identity on the vertices $l(\sigma(Y)) + 1,\ldots,j$ and zero elsewhere, and the component map $M_Q \rightarrow X$ is the identity on the vertices $i+1,\ldots,r(\sigma(X))$ and zero elsewhere. The maps through $M_{Q'}$ are similar, but with the necessary signs to make the sequence exact.

\begin{proposition}\label{prop:extension2}
        Let $\eta_2(x)$ be as in Equation~\eqref{extension2}. Then $\eta_2(x)$ is a nonsplit short exact sequence in $\Ext^1_{\Pi(A_n,S)}(X,Y)$.
\end{proposition}

\begin{proof}
    The result is again shown for $S = [2,n-1]_\mathbb{Z}$ in \cite[Theorem~8.5]{BDMTY}, and thus follows for $S$ arbitrary as in the proof of Proposition~\ref{prop:extension1}.
\end{proof}

\begin{remark} 
        In the setup of Proposition~\ref{prop:extension2}, the case where $\textnormal{left}(Y)$ is empty corresponds to the case where $l(\sigma(Y)) = i$. If we are not in this case, the condition in Lemma~\ref{lem:crossing_exists}(d) tells us there must be an arrow $i \leftarrow i+1$ in our quiver-drawing of $Y$. The situation at the other endpoint is similar.
\end{remark}

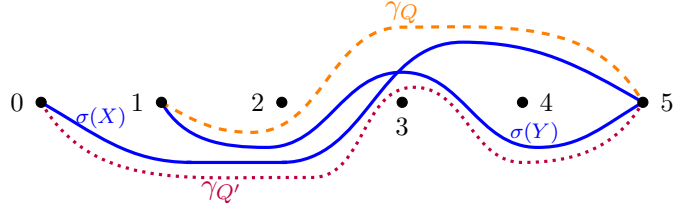
\begin{figure}
\begin{tikzpicture}[scale=0.8]

\draw[blue,very thick,smooth] (-2,0) [out = -60,in = 180] to (-0.25,-0.75) [out = 0,in = 180] to (2,0.5) [out = 0,in = 180] to (4.25,-0.75) [out = 0,in = -150] to(6,0);

\draw[blue,very thick,smooth] (-4,0) [out = -30,in = 180] to (-1.5,-1) [out = 180,in = 0] to (0,-1) [out = 0,in = 180] to (3,1) [out = 0,in = 150] to (6,0);

\draw[orange,dashed,very thick,smooth] (-2,0) [out = -30,in = 180] to (-0.5,-0.5) [out = 0,in = 180] to (2,1.25) [out = 0,in = 120] to (6,0);

\draw[purple,dotted,very thick,smooth] (-4,0) [out = -60,in = 180] to (0,-1.25) [out = 0,in = 180] to (0.5,-1.25) [out = 0,in = 180] to (2.2,0.25) [out = 0,in = 180] to (4,-1) [out = 0,in = -120] to (6,0);

\node at (2,1.5) {{\color{orange}$\gamma_{Q}$}};
\node at (-1,-1.5) {{\color{purple}$\gamma_{Q'}$}};
\node at (4.2,-0.5) {{\color{blue}\scriptsize$\sigma(Y)$}};
\node at (-3,-0.25) {{\color{blue}\scriptsize$\sigma(X)$}};

\node at (-4.4,0) {{\small 0}};
\node at (-2.4,0) {{\small 1}};
\node at (-0.4,0) {{\small 2}};
\node at (2,-0.4) {{\small 3}};
\node at (4.4,0) {{\small 4}};
\node at (6.4,0) {{\small 5}};

\filldraw  (0,0) circle (2.4pt)
 (-2,0) circle (2.4pt)
 (-4,0) circle (2.4pt)
 (2,0) circle (2.4pt)
  (4,0) circle (2.4pt)
 (6,0) circle (2.4pt)
 ;

\end{tikzpicture}
\caption{The arcs in Example~\ref{ex:extension2}. The arc $\sigma(X)$ has endpoints 0 and 5 and the arc $\sigma(Y)$ has endpoints 1 and 5.}\label{fig:extension2}
\end{figure}

\begin{example}\label{ex:extension2}
    Let $X = \text{\tiny $\begin{matrix}1 \ 5\\ 2\ 4\\3\end{matrix}$}$ and $Y = \text{\tiny $\begin{matrix}2 \ 4\\ \ \ 3 \ 5\end{matrix}$}$ in $\brick(\Pi(A_5,S))$. The corresponding arcs are shown in Figure~\ref{fig:extension2}. There is a unique nontrivial crossing, which is directed from $\sigma(Y)$ to $\sigma(X)$. The corresponding nodes $i$ and $j$ are $i = 1 = l(\sigma(Y))$ and $j = 4$. The modules $M_Q$ and $M_{Q'}$ in Equation~\eqref{extension2} are then 
    given by $M_Q = \text{\tiny $\begin{matrix}1 \ \ \ \ \ \\2 \ 4\\ \ \ 3 \ 5\end{matrix}$}$ and $M_{Q'} = \text{\tiny $\begin{matrix}\ \ \ \ 5\\ 2\ 4\\3\end{matrix}$}.$
    The corresponding arcs $\gamma_Q$ and $\gamma_{Q'}$ are shown in Figure~\ref{fig:extension2}. Note that these arcs can be formed as follows. Cut $\sigma(X)$ and $\sigma(Y)$ at their point of intersection. Glue the left half of $\sigma(Y)$ with the right half of $\sigma(X)$ and perturb upwards to obtain $\gamma_Q$. Similarly, glue the left of $\sigma(X)$ with the right half of $\sigma(Y)$ and perturb downward to obtain $\gamma_{Q'}$.
\end{example}

Finally, suppose there is a nontrivial crossing $x$ directed from $\sigma(X)$ to $\sigma(Y)$. Let $i < j$ be the unique integers satisfying the properties (a)-(e) in Lemma~\ref{lem:crossing_exists} for $\gamma_1 = \sigma(X)$ and $\gamma_2 = \sigma(Y)$. Denote by $Q$ and $Q'$ the quivers corresponding to $X$ and $Y$, respectively. For each node $i \in [n]$, if they exist then we denote by $v_i$ and $v'_i$ the vertices labeled by $i$ in $Q$ and $Q'$, respectively. If $Q$ has a vertex $v_i$ labeled by $i$, draw an arrow $v_i \rightarrow v'_{i+1}$. Otherwise, $Q'$ has a vertex $v'_i$ labeled by $i$, so draw an arrow $v_{i+1} \rightarrow v'_i$. Similarly, if $Q$ has a vertex $v_{j+1}$ labeled by $j+1$, draw an arrow $v_{j+1} \rightarrow v'_j$. Otherwise, $Q'$ has a vertex $v'_{j+1}$ labeled by $j+1$, so draw an arrow $v_j \rightarrow v'_{j+1}$. Moreover, for all nodes $i < k < j$, if $Q$ contains an arrow $v_k \rightarrow v_{k+1}$, then draw an arrow $v_{k+1}\rightarrow v'_k$. Similarly, if $Q$ contains an arrow $v_k \leftarrow v_{k+1}$, draw an arrow $v_k \rightarrow v'_{k+1}$. The resulting quiver $Q''$ resembles:
    $$
    \begin{tikzcd}[column sep = 0.35cm]
        \textnormal{left}(X)\arrow[dr,dashed] & v_{i+1} \arrow[l,dashed] \arrow[r] & v_{i+2}\arrow[dl] \arrow[r] & v_{i+3}\arrow[dr]\arrow[dl] & \arrow[l] \cdots \arrow[r] & v_{j-2}\arrow[dl]\arrow[r] & v_{j-1}\arrow[dl]\arrow[dr] & v_j \arrow[l] \arrow[r,dashed] & \textnormal{right}(X)\arrow[dl,dashed]\\
        \textnormal{left}(Y) \arrow[r,dashed] & v'_{i+1}  \arrow[r] & v'_{i+2} \arrow[r] & v'_{i+3} & \arrow[l] \cdots \arrow[r] & v'_{j-2}\arrow[r] & v'_{j-1} & v'_j \arrow[l] & \textnormal{right}(Y), \arrow[l,dashed]
    \end{tikzcd}
    $$
    where if $\textnormal{left}(X)$ is empty we replace the arrow from $v_i$ to $v'_{i+1}$ with one from $v_{i+1}$ to $v'_i$, and similarly on the right of the diagram. (See also Example~\ref{ex:extension3} for a specific example.) Now let
    \begin{equation}\label{extension3}
        \eta_3(x) = \qquad 0 \rightarrow Y \rightarrow M_{Q''} \rightarrow X \rightarrow 0,
    \end{equation}
    where the first map is given by the inclusion of $Q'$ as the bottom row of $Q''$ and the second is given by the projection from the top row of $Q''$ to $Q$. 

\begin{proposition}\label{prop:extension3}
    Let $M_{Q''}$ and $\eta_3(x)$ be as in Equation~\eqref{extension3} above. Then $M_{Q''} \in \mods \Pi(A_n,S)$ if and only if $S \cap [i+1,j]_\mathbb{Z} = \emptyset$. Moreover, if $M_{Q''} \in \mods \Pi(A_n,S)$, then $\eta_3(x)$ is a nonsplit short exact sequence in $\Ext^1_{\Pi(A_n,S)}(X,Y)$.
\end{proposition}

\begin{proof}
    By the assumptions (a)-(e) in Lemma~\ref{lem:crossing_exists}, we observe the following:
    \begin{itemize}
        \item $\textnormal{left}(X)$ (resp.\ $\textnormal{left}(Y)$) is empty if and only if $l(\sigma(X)) = i$ (resp.\ $l(\sigma(Y)) = i$). In particular, at least one of these is nonempty. The situation is similar for $\textnormal{right}(X)$ and $\textnormal{right}(Y)$.
        \item If $\textnormal{left}(X)$ (resp.\ $\textnormal{right}(Y)$) is not empty, then $Q$ contains an arrow $i \leftarrow i+1$ (resp.\ $Q'$ contains an arrow $i' \rightarrow (i+1)'$). The situation is similar for $\textnormal{right}(X)$ and $\textnormal{right}(Y)$. 
        \item If $k \in (i,j)\cap \mathbb{Z}$, then $\sigma(X)$ and $\sigma(Y)$ must pass on the same side of $k$. In particular, the arrows connecting $k$ with $k+1$ in $Q$ and $k'$ with $(k+1)'$ in $Q'$ must point in the same direction.
    \end{itemize}
    This shows that $M_{Q''}$ satisfies the relations defining $\Pi(A_n)$, and more precisely that $M_{Q''}\in \mods \Pi(A_n,S)$ if and only if $S \cap [i+1,j]_\mathbb{Z} = \emptyset$. The fact that $\eta_3(x)$ is a short exact sequence then follows immediately. Finally, the fact that this sequence is not split follows from the fact that $\eta_3(x) \notin \Ext^1_{\Pi(A_n,[2,n-1]_\mathbb{Z})}(X,Y)$, which contains the trivial extension.
\end{proof}

\begin{example}\label{ex:extension3}
    Let $X = \text{\tiny $\begin{matrix}2 \ 4\\ \ \ 3 \ 5\end{matrix}$}$ and $Y = \text{\tiny $\begin{matrix}1 \ 5\\ 2\ 4\\3\end{matrix}$}$ be  in $\brick(\Pi(A_5,S))$. (Note that this is the opposite assignment of $X$ and $Y$ from Example~\ref{ex:extension2}.) As in Example~\ref{ex:extension2}, there is a unique nontrivial crossing, directed from $\sigma(X)$ to $\sigma(Y)$, with corresponding nodes $i = 1 = l(\sigma(X))$ and $j = 4$. The $M_{Q''}$ in Equation~\eqref{extension3} then has
    $$Q'' = 
    \begin{tikzcd}
        & 2\arrow[r]\arrow[dl] & 3\arrow[dl]\arrow[dr] & 4\arrow[l]\arrow[r] & 5\arrow[dl]\\
        1 \arrow[r] & 2 \arrow[r] & 3 & 4\arrow[l] & 5.\arrow[l]
    \end{tikzcd}
    $$
    (In particular, $\textnormal{left}(X)$ is empty in this case.) We see that $M_{Q''}$ is a module over $\Pi(A_5) = \Pi(A_5,\emptyset)$, but not over $\Pi(A_5,\{2,3,4\})$. Moreover, we have that $Y$ is injective over $\Pi(A_5,\{2,3,4\})$, but not over $\Pi(A_5)$. Thus the short exact sequence in equation \ref{extension3} could not possibly exist over $\Pi(A_5,\{2,3,4\})$. Finally, we note that the module $M_{Q''}$ is not indecomposable in this case: take the basis $\{(1_K,0), (1_K,-1_K)\}$ for $M_{Q''}(5)$, where the first coordinate represents the projection to the top copy of ``5'' and the second coordinate the projection to the bottom copy of ``5''. Then the subspace spanned by $(1_K,-1_K)$ constitutes a copy of the simple module $S(5)$ which is both a submodule and a quotient module of $M_{Q''}$. There are, however, examples where the middle term of this extension is indecomposable. Indeed, this is the case if we remove the final composition factor ``5'' from~$X$.
\end{example}

We now prove our main result of this section.

\begin{theorem}\label{thm:ext_space}
    Let $X, Y \in \brick(\Pi(A_n))$.
    \begin{enumerate}
        \item Suppose that $\sigma(X)$ and $\sigma(Y)$ have a contested endpoint $x$. Then $\Ext^1_{\Pi(A_n,[2,n-1]_\mathbb{Z})}(X,Y) = \Ext^1_{\Pi(A_n)}(X,Y)$ is 1-dimensional and spanned by the short exact sequence $\eta_1(x)$ from Equation~\eqref{extension1}.
        \item Suppose that $\sigma(X)$ and $\sigma(Y)$ do not have a contested endpoint. Choose representatives of $\sigma(X)$ and $\sigma(Y)$ which intersect a minimal number of times, let $\mathcal{X}$ be the set of nontrivial crossings between $\sigma(X)$ and $\sigma(Y)$, and let $\mathcal{X}_- \subseteq \mathcal{X}$ be the subset of crossings which are directed from $\sigma(Y)$ to $\sigma(X)$. Then:
        \begin{enumerate}
            \item $\{\eta_2(x) \mid x \in \mathcal{X}_-\}$ is a basis of $\Ext^1_{\Pi(A_n,[2,n-1]_\mathbb{Z})}(X,Y)$, where each $\eta_2(x)$ is as defined in Equation~\eqref{extension2}.
            \item $\{\eta_2(x) \mid x \in \mathcal{X}_-\} \cup \{\eta_3(x) \mid x \in \mathcal{X} \setminus \mathcal{X}_-\}$ is a basis of $\Ext^1_{\Pi(A_n)}(X,Y)$, where each $\eta_2(x)$ is as defined in Equation~\eqref{extension2} and each $\eta_3(x)$ is defined as in Equation~\eqref{extension3}.
        \end{enumerate}
    \end{enumerate}
\end{theorem}

\begin{proof}
    (1) This is contained in Proposition~\ref{prop:extension1}.

    (2a) This result is essentially contained in \cite[Section~6.2]{BaH}, where it is deduced from \cite[Theorem~8.5]{BDMTY} and \cite[Theorem~A]{CPS}. Indeed, in \cite{BDMTY,CPS}, the authors give explicit bases for the Hom- and Ext-spaces between indecomposable ``string modules'' over gentle algebras. (The basis for the Hom-space is actually much more classical. It was first established in higher generality in \cite{CB_string} and then specialized to gentle algebras in \cite{schroer}.) In the terminology of \cite{BDMTY}, the basis elements of $\Hom_{\Pi(A_n,[2,n-1]_\mathbb{Z})}(Y,X)$ are called ``graph maps''. The image of a graph map $f: Y \rightarrow X$ is again a ``string module'', and is therefore indecomposable. Moreover, since $Y$ has dimension at most 1 at every vertex, $\mathrm{im}(f)$ is actually a brick. Now the graph map $f$ is called ``two-sided'' if there exists $i < j < k \in [n] \setminus\{0\}$ such that (a) $X(i) \oplus Y(i) \neq 0$ and $\mathrm{im}(f)(i) = 0$, (b) $\mathrm{im}(f)(j) \neq 0$, and (c) $X(k) \oplus Y(k) \neq 0$ and $\mathrm{im}(f)(k) = 0$. As is implicit in \cite[Theorem~6.2.5]{BaH}, the two-sided graph maps from $Y$ to $X$ are those corresponding (in the proof of Proposition~\ref{prop:hom_space}) to the nontrivial crossings directed from $\sigma(Y)$ to $\sigma(X)$. It then follows from \cite[Theorem~6.5]{BDMTY} or \cite[Theorem~A]{CPS} that the nonsplit sequences $\{\eta_2(x) \mid x \in \mathcal{X}_-\}$ form a basis of $\Ext^1_{\Pi(A_n,[2,n-1]_\mathbb{Z})}(X,Y)$.

    (2b) Denote $\mathcal{B}_+ = \{\eta_3(x) \mid x \in \mathcal{X} \setminus \mathcal{X}_-\}$ and $\mathcal{B}_- = \{\eta_2(x) \mid x \in \mathcal{X}_-\}$, and note that $\mathcal{B}_- \cap \mathcal{B}_+ = \emptyset$. By Corollary~\ref{cor:ext}, we have that $\left|\mathcal{X}\right| = \left|\mathcal{B}_+ \cup \mathcal{B}_-\right| = \dim\Ext^1_{\Pi(A_n)}(X,Y)$. Thus we need only show that the exact sequences in $\mathcal{B}_+ \cup \mathcal{B}_-$ are linearly independent.

    Let $\mathcal{C} \subseteq \mathcal{B}_+$. We will show that the vectors in $\mathcal{C} \cup \mathcal{B}_-$ are linearly independent by induction on $|\mathcal{C}|$. For the base case, we note that the vectors in $\mathcal{B}_-$ form a basis of $\Ext^1_{\Pi(A_n,[2,n-1]_\mathbb{Z})}(X,Y) \subseteq \Ext^1_{\Pi(A_n)}(X,Y)$ by (2a). 
    
    Now suppose that $|\mathcal{C}| = k > 1$, and that the result holds for $k-1$. Choose some $\eta_3(x) \in \mathcal{C}$, and let $i < j$ be the nodes satisfying conditions (a-e) in Lemma~\ref{lem:crossing_exists} for the crossing $x$ (with $\gamma_1 = \sigma(X)$ and $\gamma_2 = \sigma(Y)$). Then Lemma~\ref{lem:crossing_exists} implies that $x$ is the only intersection point between $\sigma(X)$ and $\sigma(Y)$ in the interval $(i,j)$. Thus we have $\eta_3(x) \notin \Ext^1_{\Pi(A_n,[i+1,j]_\mathbb{Z})}(X,Y)$ and that $\mathcal{B}_- \cup (\mathcal{C} \setminus \{\eta_3(x)\}) \subseteq \Ext^1_{\Pi(A_n,[i+1,j]_\mathbb{Z}))}(X,Y)$. The linear independence of $\mathcal{B}_- \cup \mathcal{C}$ thus follows from the induction hypothesis.
\end{proof}

As an analog to Corollary~\ref{cor:hom}, Theorem~\ref{thm:ext_space} implies the following.

\begin{corollary}\label{cor:ext}
    Let $X, Y \in \brick(\Pi(A_n,[2,n-1]_\mathbb{Z}))$. Then $$\left|\dim\Ext^1_{\Pi(A_n,[2,n-1]_\mathbb{Z})}(X,Y) - \dim\Ext^1_{\Pi(A_n,[2,n-1]_\mathbb{Z})}(Y,X)\right| \leq 1.$$
\end{corollary}

\begin{proof}
    Choose representatives of the arcs $\sigma(X)$ and $\sigma(Y)$ which intersect a minimal number of times. If $\sigma(X)$ and $\sigma(Y)$ have a contested endpoint, then $\dim\Ext^1_{\Pi(A_n,[2,n-1]_\mathbb{Z})}(X,Y) = 1 = \dim\Ext^1_{\Pi(A_n,[2,n-1]_\mathbb{Z})}(Y,X)$ by Theorem~\ref{thm:ext_space}. Otherwise, let $x_1 < \cdots < x_k$ be the nontrivial crossings between $\sigma(X)$ and $\sigma(Y)$. Then for $1 \leq i < k$, the minimality of the representatives implies that if $x_i$ is directed from $\sigma(X)$ to $\sigma(Y)$, then $x_{i+1}$ is directed from $\sigma(Y)$ to $\sigma(X)$, and vice versa. The result thus follows from Theorem~\ref{thm:ext_space}.
\end{proof}

We conclude this section by tabulating three additional consequences of Theorem~\ref{thm:ext_space}. The first will be critical in the next two sections, and can also be found in \cite[Theorem~1.2]{IRRT} and \cite[Proposition~3.7]{mizuno2}

\begin{corollary}\label{cor:bricks_K}
    Let $S \subseteq [2,n-1]_\mathbb{Z}$ and let $X \in \brick(\Pi(A_n,S))$. Then $\mathrm{End}_{\Pi(A_n,S)}(X) \cong K$ and $\Ext^1_{\Pi(A_n,S)}(X,X) = 0$.
\end{corollary}

\begin{proof}
    We can choose two representatives $\gamma$ and $\gamma'$ for the arc $\sigma(X)$ so that $\gamma$ and $\gamma'$ share both of their endpoints and do not otherwise intersect. As Hom-crossings, one of the endpoints will be directed from $\gamma$ to $\gamma'$ and the other from $\gamma'$ to $\gamma$. The result then follows immediately from Proposition~\ref{prop:hom_space} and Theorem~\ref{thm:ext_space}(2).
\end{proof}

Now recall that a module $M \in \mods \Pi(A_n,S)$ is called \emph{$\tau$-rigid} if $\Hom_{\Pi(A_n,S)}(M,\tau_S M) = 0$, where we use $\tau_S$ as shorthand for the Auslander-Reiten translation in $\mods \Pi(A_n,S)$. Such modules are critical in the study of \emph{$\tau$-tilting theory} \cite{AIR}, and will be used in our sequel work describing the \emph{$\tau$-exceptional sequences} over $\Pi(A_n)$. By a well-known result of \cite{AS}, we have that $\Hom_{\Pi(A_n,S)}(M,\tau_S N) = 0$ if and only if $\Ext^1_{\Pi(A_n,S)}(N,M') = 0$ for every indecomposable quotient $M'$ of $M$. Thus we can prove the following.

\begin{corollary}\label{cor:tau_rigid}
    Let $X \in \brick(\Pi(A_n))$ and write $\wrd(\sigma(X)) = \left(l_X,s_{l_X+1}\cdots s_{r_X}\right)$. Then:
    \begin{enumerate}
        \item $X$ is $\tau$-rigid in $\mods\Pi(A_n,[2,n-1]_\mathbb{Z})$.
        \item The following are equivalent:
        \begin{enumerate}
            \item $X$ is $\tau$-rigid in $\mods \Pi(A_n)$.
            \item $X$ has a simple socle.
            \item There exists $i \in [l_X+1,r_X]_\mathbb{Z}$ such that $s_j = \textnormal{u}$ for all $l_X < j < i$ and $s_j = \textnormal{o}$ for all $i \leq j < r_X$.
        \end{enumerate}
    \end{enumerate}
\end{corollary}

\begin{proof}
    (1) Let $Y$ be an indecomposable quotient of $X$ in $\mods \Pi(A_n)$, and suppose that $Y 
 \neq X$. Propositions~\ref{prop:quotients} and~\ref{prop:hom_space} then imply that (i) $Y \in \mods\Pi(A_n,[2,n-1]_\mathbb{Z})$, (ii) $\sigma(Y)$ is a quotient arc of $\sigma(X)$, and (iii) there is a unique Hom-crossing between $\sigma(X)$ and $\sigma(Y)$, and this Hom-crossing is directed from $\sigma(X)$ to $\sigma(Y)$. It then follows from Theorem~\ref{thm:ext_space}(2a) that $\Ext^1_{\Pi(A_n,[2,n-1]_\mathbb{Z}}(X,Y) = 0$. The result thus follows from the fact that $\Ext^1_{\Pi(A_n,[2,n-1]_\mathbb{Z}}(X,X) = 0$ established in Corollary~\ref{cor:bricks_K}.

    (2) By Theorem~\ref{thm:ext_space}(2b), the proof of (1) implies that $X$ is $\tau$-rigid in $\mods \Pi(A_n)$ if and only if every quotient arc of $\sigma(X)$ has a shared endpoint with $\sigma(X)$.

    $(a \implies b)$: Suppose that $\mathrm{soc}(X)$ is not simple. Since $X$ has dimension at most 1 at every vertex, this means there exist $i < j \in [n] \setminus \{0\}$ such that both $S(i)$ and $S(j)$ are submodules of $X$. Now recall that $\wrd(\sigma(S(i))) = (i-1,\textnormal{e})$, and likewise for $S(j)$. Thus by Proposition~\ref{prop:quotients}, we have that $s_{i} = \textnormal{o}$ and $s_{j-1} = \textnormal{u}$. In particular, we can choose $k \in [i+1,j-1]_\mathbb{Z}$ minimal such that $s_{k-1} = \textnormal{o}$ and $s_k = \textnormal{u}$. Then the arc $\sigma(S(k))$ is a quotient arc of $\sigma(X)$ which does not have a shared endpoint with $\sigma(X)$. We conclude that $X$ is not $\tau$-rigid in $\mods \Pi(A_n)$.
    
    $(b \implies c)$: Suppose $X$ has a simple socle. By Proposition~\ref{prop:quotients}, this means there is a unique $i \in [l_X+1,r_X]_\mathbb{Z}$ such that $\sigma(S(i))$ is a submodule arc of $\sigma(X)$. In particular, $s_i \neq \textnormal{u}$ and either $i - 1 = l_X$ or $s_{i-1} = \textnormal{u}$. Since $i$ is unique, it follows that $s_j = \textnormal{u}$ for all $l_X < j < i$ and that $s_j = \textnormal{o}$ for all $i \leq j < r_X$.

    $(c \implies a)$: Suppose $X$ is not $\tau$-rigid. Then there exists a quotient arc $\gamma$ of $\sigma(X)$ which does not have a shared endpoint with $\sigma(X)$. Thus $s_{l(\gamma)} = \textnormal{o}$ and $s_{r(\gamma)} = \textnormal{u}$. We conclude that condition (c) must be false.
\end{proof}

We note that Corollary~\ref{cor:tau_rigid}(1) is also a consequence of the ``brick-$\tau$-rigid correspondence'' of \cite{DIJ}. Indeed, this result implies that, if an algebra has finitely many bricks, then its bricks and its indecomposable $\tau$-rigid modules are in bijection with one another. Since every indecomposable module over $\Pi(A_n,[2,n-1]_\mathbb{Z})$ is a brick, (1) then follows immediately.

We conclude this section with a brief discussion of the ``almost rigid modules'' of \cite{BCSGS,BGMS}.

\begin{remark}\label{rem:MAR}
    In \cite[Definition~6.1]{BGMS}, a basic module $M$ over a path algebra of type A is said to be \emph{almost rigid} if for any two indecomposable direct summands $M_1, M_2$ of $M$ both (i) any nonsplit short exact sequence of the form $0 \rightarrow M_2 \rightarrow E \rightarrow M_1 \rightarrow 0$ has $E$ indecomposable, and (ii) $\dim \Ext^1(M_1,M_2) \leq 1$. If in addition $M$ is not a direct summand of a larger almost rigid module then $M$ is said to be \emph{maximal almost rigid}. These definitions are being extended to all gentle algebras in the forthcoming work \cite{BCSGS}.

    For $X \neq Y \in \brick(\Pi(A_n,[2,n-1]_\mathbb{Z}))$, it follows immediately from Theorem~\ref{thm:ext_space} that $X \oplus Y$ is almost rigid if and only if the corresponding arcs $\sigma(X)$ and $\sigma(Y)$ have no nontrivial crossings. By taking ``noncrossing'' to mean ``without nontrivial crossings'', this yields a bijection between the (maximal) almost rigid modules over $\Pi(A_n,[2n-1])$ and (maximal) collections of pairwise noncrossing arcs in $\arc(n)$.
\end{remark}


\section{Weak exceptional sequences}\label{sec:weak}

In this section, we classify the weak exceptional sequences over the algebras $\Pi(A_n,S)$ using arcs. In particular, we show that the set of weak exceptional sequences does not depend on the choice of the subset $S \subseteq [2,n-1]_\mathbb{Z}$. We then use our characterization to show that if $n \neq 1$ then the maximum length of a weak exceptional sequence over $\Pi(A_n,S)$ is $2n-2$, see Corollary~\ref{cor:planar2}. We first observe the following, which is also shown in \cite[Section~6.2]{BaH} (see also Remark~\ref{rem:weak}(2)).

\begin{lemma}\label{lem:same_weak}
    Let $(Y,X)$ be a sequence of bricks in $\brick(\Pi(A_n))$. Choose representatives of the arcs $\sigma(X)$ and $\sigma(Y)$ which intersect a minimal number of times. Then the following are equivalent.
    \begin{enumerate}
        \item There exists $S \subseteq [2,n-1]_\mathbb{Z}$ such that $(Y,X)$ is a weak exceptional pair in $\mods \Pi(A_n,S)$.
        \item $(Y,X)$ is a weak exceptional pair in $\mods \Pi(A_n,S)$ for all $S \subseteq [2,n-1]_\mathbb{Z}$.
        \item Either (a) $\sigma(X)$ and $\sigma(Y)$ do not intersect or (b) $\sigma(X)$ and $\sigma(Y)$ have one shared endpoint, as a Hom-crossing this endpoint is directed from $\sigma(Y)$ to $\sigma(X)$, and $\sigma(X)$ and $\sigma(Y)$ do not otherwise intersect.
    \end{enumerate}
\end{lemma}

\begin{proof}
    $(1 \implies 3)$: Suppose $(Y,X)$ is a weak exceptional pair in $\mods \Pi(A_n,S)$. Then by Proposition~\ref{prop:hom_space} there are no Hom-crossings directed from $\sigma(X)$ to $\sigma(Y)$. By Corollary~\ref{cor:hom}, this means there is at most one Hom-crossing directed from $\sigma(Y)$ to $\sigma(X)$. Moreover, since $\Ext^1_{\Pi(A_n,[2,n-1]_\mathbb{Z})}(X,Y) \subseteq \Ext^1_{\Pi(A_n,S)}(X,Y)$, this Hom-crossing cannot be a nontrivial crossing by Theorem~\ref{thm:ext_space}. We conclude that (3a) and (3b) are the only possibilities.

    $(3 \implies 2)$: Suppose (3) holds. It follows from Proposition~\ref{prop:hom_space} and Theorem~\ref{thm:ext_space} that $\Hom_{\Pi(A_n)}(X,Y) = 0 = \Ext^1_{\Pi(A_n)}(X,Y)$. Then for any $S$, the fact that $\Ext^1_{\Pi(A_n,S)}(X,Y) \subseteq \Ext^1_{\Pi(A_n)}(X,Y)$ implies that $(Y,X)$ is a weak exceptional pair over $\Pi(A_n,S)$.

    $(2 \implies 1)$: trivial.
\end{proof}

An immediate consequence is the following.

\begin{corollary}\label{cor:same_weak}
    The weak exceptional sequences in $\mods \Pi(A_n,S)$ do not depend on the choice of the subset $S \subseteq [2,n-1]_\mathbb{Z}$.
\end{corollary}

In order to establish a combinatorial model for these weak exceptional sequences, we consider the following definitions.

\begin{definition}\label{def:clockwise}
    Let $\gamma_1, \gamma_2 \in \arc(n)$. We say that $\gamma_2$ is \emph{clockwise} of $\gamma_1$ if all of the following hold.
    \begin{enumerate}
        \item The arcs $\gamma_1$ and $\gamma_2$ have one shared endpoint and zero Ext-crossings.
        \item The shared endpoint between $\gamma_1$ and $\gamma_2$ is directed from $\gamma_2$ to $\gamma_1$ as a Hom-crossing.
    \end{enumerate}
\end{definition}

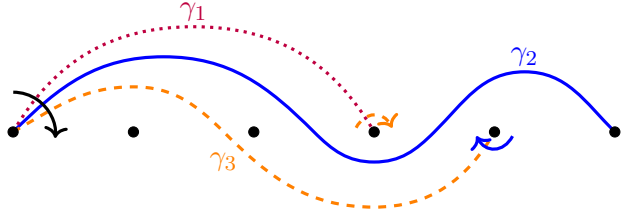
\begin{figure}
\begin{tikzpicture}[scale=0.8]

\draw [blue, very thick] 
    (-2,0) [out = 45,in = 180] to (0.5,1.25) [out = 0,in = 135] to (3,0) [out = -45,in = 180] to (4,-0.5) [out = 0,in = -135] to (5,0) [out = 45,in = 180] to (6.5,1) [out = 0,in = 135] to (8,0);

\draw[orange, very thick,dashed]
    (-2,0) [out = 30,in = 180] to (-0,0.75) [out = 0,in = 135] to (1.5,0) [out = -45,in = 180] to (4,-1.25) [out = 0,in = -120] to (6,0);
;

\draw[purple, very thick,smooth,dotted]
    (-2,0) [out = 60,in=180] to (1,1.75) [out=0,in=120] to (4,0);

\node [purple] at (1,2) {$\gamma_1$};
\node [blue] at (6.5,1.25) {$\gamma_2$};
\node [orange] at (1.5,-0.5) {$\gamma_3$};

\draw [orange, very thick, ->,dashed] 
(4-0.3, 0.22-0.15) .. controls (4-0.15, 0.5-0.15)and(4+0.15, 0.5-0.15) .. (4+0.3, 0.22-0.15);

\draw [blue, very thick, <-] 
(6-0.3, -0.22+0.15) .. controls (6-0.15, -0.5+0.15)and(6+0.15, -0.5+0.15) .. (6+0.3, -0.22+0.15);

\filldraw  (0,0) circle (2.4pt)
 (-2,0) circle (2.4pt)
 (2,0) circle (2.4pt)
 (4,0) circle (2.4pt) 
 (6,0) circle (2.4pt)
 (8,0) circle (2.4pt)
 ;

\draw[black, very thick, <-] (-1.3, -0.1) arc (-5:90:7mm);
\end{tikzpicture}
\caption{The arc $\gamma_3$ (dashed orange) is clockwise of both $\gamma_1$ (dotted purple) and $\gamma_2$ (solid blue). At the right endpoint of $\gamma_1$, the small clockwise arrow of $\gamma_3$ points ``away from'' $\gamma_1$, and at the right endpont of $\gamma_3$, the small clockwise arrow of $\gamma_2$ points ``towards'' $\gamma_3$.}\label{fig:clockwise}
\end{figure}

\begin{remark}\label{rem:clockwise}
     Let $\gamma_1, \gamma_2 \in \arc(n)$. For $t \in \{1,2\}$, denoted $\wrd(\gamma_t) = \left(l(\gamma_t),s_{l(\gamma_t)+1}^t\cdots s_{r(\gamma_t)}\right).$
    \begin{enumerate}
        \item Suppose that $l(\gamma_1) = l(\gamma_2)$ and that $\gamma_1$ and $\gamma_2$ do not otherwise intersect. Then $\gamma_2$ is clockwise of $\gamma_1$ if and only if either (a) $r(\gamma_1) < r(\gamma_2)$ and $s_{r(\gamma_1)}^2 = \text{u}$, or (b) $r(\gamma_1) > r(\gamma_2)$ and $s_{r(\gamma_2)}^1 = \text{o}$. Visually, in case (a) the small clockwise arrow of $\gamma_1$ at $r(\gamma_2)$ points ``towards'' $\gamma_2$. Similarly, in case (b) the small clockwise arrow of $\gamma_2$ at $r(\gamma_1)$ points ``away from'' $\gamma_1$. See Figure~\ref{fig:clockwise} for an example. 
        \item Suppose that $r(\gamma_1) = r(\gamma_2)$ and that $\gamma_1$ and $\gamma_2$ do not otherwise intersect. Then $\gamma_2$ is clockwise of $\gamma_1$ if either (a) $l(\gamma_1) < l(\gamma_2)$ and $s_{l(\gamma_2)}^1 = \text{o}$, or (b) $l(\gamma_1) > l(\gamma_2)$ and $s_{l(\gamma_1)}^2 = \text{u}$. These leads to a similar visualization as in the $l(\gamma_1) = l(\gamma_2)$ case.
        \item If $\gamma_1$ and $\gamma_2$ have a contested endpoint, then $\gamma_2$ is not clockwise of $\gamma_1$ and vice versa.
    \end{enumerate}
\end{remark}

\begin{definition}\label{def:cw_arc}
    Let $\omega 
 = (\gamma_k,\ldots,\gamma_1)$ be an ordered set of arcs in $\arc(n)$. We say that $\omega$ is a \emph{clockwise-ordered arc diagram (}of length $k$) if the following hold for all $1 \leq i < j \leq k$.
    \begin{enumerate}
        \item There exist representatives of $\gamma_i$ and $\gamma_j$ which have at most one shared endpoint and do not otherwise intersect.
        \item If $\gamma_i$ and $\gamma_j$ have a shared endpoint then $\gamma_j$ is clockwise of $\gamma_i$.
    \end{enumerate}
    We denote by $\arc_{cw}(n,k)$ the set of combinatorial equivalence classes of clockwise-ordered arc diagrams on $n+1$ of length $k$ and $\arc_{cw}(n) := \bigcup_{k \in \mathbb{N}}\arc_{cw}(n,k)$.
\end{definition}

We note that it does not follow \emph{a priori} that there exist \emph{fixed} representatives of the arcs in a given clockwise-ordered arc diagram which intersect only at shared endpoints. Indeed, we have only assumed such representatives can be chosen in a pairwise fashion. It is, however, true that such representatives always exists, as the following shows.
 \begin{proposition}\label{prop:planar}
    Let $\omega = (\gamma_k,\ldots,\gamma_1) \in \arc_{cw}(n)$ be a clockwise-ordered arc diagram. Then the arcs in $\omega$ together form a simple planar bipartite graph on $n+1$ nodes.
\end{proposition}

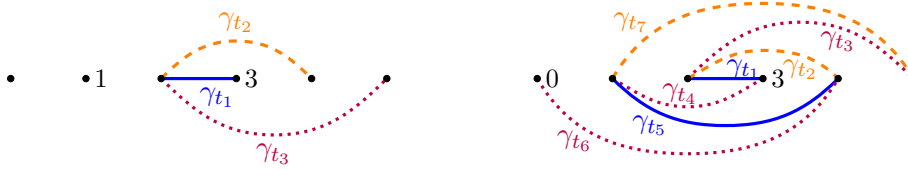
\begin{figure}
\begin{tikzpicture}[scale=0.5]

\draw [blue, very thick,smooth] 
(2,0) -- (4,0);

\draw [orange, very thick,smooth,dashed] 
(2,0) [out = 45, in = 180] to (4,1) [out = 0, in = 135] to (6,0);

\draw [purple, very thick,smooth,dotted] 
(2,0) [out = -45, in = 180] to (5,-1.5) [out = 0, in = -135] to (8,0);

\node[color = blue] at (3.5,-0.5){$\gamma_{t_1}$};

\node[color = purple] at (5,-2){$\gamma_{t_3}$};

\node[color = orange] at (4,1.5){$\gamma_{t_2}$};
        
\filldraw  (0,0) circle (2.4pt)
 (-2,0) circle (2.4pt)
 (2,0) circle (2.4pt)
 (4,0) circle (2.4pt) 
 (6,0) circle (2.4pt)
 (8,0) circle (2.4pt)
 ;

\node at (0.4,0) {{\small 1}};
\node at (4.4,0) {{\small 3}};

\begin{scope}[shift = {(14,0)}]
\draw [blue, very thick,smooth] 
(2,0) -- (4,0);

\draw [orange, very thick,smooth,dashed] 
(2,0) [out = 30, in = 180] to (4,0.75) [out = 0, in = 135] to (6,0);

\draw [purple, very thick,smooth,dotted] 
(2,0) [out = 45, in = 180] to (5,1.5) [out = 0, in = 135] to (8,0);

\draw [purple, very thick,smooth,dotted] 
(0,0) [out = -30, in = 180] to (2,-0.75) [out = 0, in = -150] to (4,0);

\draw [blue, very thick,smooth] 
(0,0) [out = -45, in = 180] to (3,-1.25) [out = 0, in = -135] to (6,0);

\draw [purple, very thick,smooth,dotted] 
(-2,0) [out = -60, in = 180] to (2,-2) [out = 0, in = -120] to (6,0);

\draw [orange, very thick,smooth,dashed] 
(0,0) [out = 60, in = 180] to (4,2) [out = 0, in = 120] to (8,0);

\node[color = blue] at (3.5,0.25){$\gamma_{t_1}$};

\node[color = purple] at (6,1){$\gamma_{t_3}$};

\node[color = orange] at (5,0.25){$\gamma_{t_2}$};

\node[color = purple] at (1.8,-0.4){$\gamma_{t_4}$};

\node[color = blue] at (1,-1.25){$\gamma_{t_5}$};

\node[color = purple] at (-1,-1.65){$\gamma_{t_6}$};

\node[color = orange] at (0.5,1.5){$\gamma_{t_7}$};
        
\filldraw  (0,0) circle (2.4pt)
 (-2,0) circle (2.4pt)
 (2,0) circle (2.4pt)
 (4,0) circle (2.4pt) 
 (6,0) circle (2.4pt)
 (8,0) circle (2.4pt)
 ;

\node at (-1.6,0) {{\small 0}};
\node at (4.4,0) {{\small 3}};
\end{scope}

\end{tikzpicture}
\caption{An illustration of the proof of Proposition~\ref{prop:planar}. In the left (resp.\ right) diagram, it is impossible to connect the nodes 1 and 3 (resp.\ 0 and 3) with an arc without introducing a nontrivial crossing.}\label{fig:planar}
\end{figure}

\begin{proof}
    Let $G(\omega)$ denote the corresponding graph. First note that this graph is simple by the assumption that any pair of arcs have at most one shared endpoint. Moreover, since $\omega \in \arc_{cw}(n)$, there cannot exist a node which is both the left endpoint of one arc in $\omega$ and the right endpoint of another (see Remark~\ref{rem:clockwise}(3)). This implies that $G(\omega)$ is bipartite. To show that $G(\omega)$ is planar, it therefore suffices to show that $G(\omega)$ does not contain the complete bipartite graph $K_{3,3}$ as a subgraph.

    Assume for a contradiction that $G(\omega)$ contains a copy of $K_{3,3}$. Then by deleting nodes and shifting indices as necessary, we can assume that $n = 5$ and $G(\omega) = K_{3,3}$ with bipartition $\{\{0,1,2\},\{3,4,5\}\}$. Up to vertical symmetry, this means that $\omega$ in particular contains arcs $\gamma_{t_1}$ and $\gamma_{t_2}$ characterized by $\wrd(\gamma_{t_1}) = (2,\text{e}) =: \left(l^1,s_3^1\right)$ and $\wrd(\gamma_{t_2}) = (2,\text{oe}) =: \left(l^2,s_3^2s_4^2\right)$. There are then three cases to consider.

    Suppose first that $\omega$ contains an arc $\gamma_{t_3}$ with $\wrd(\gamma_{t_3}) = (2,\text{uue}) =: \left(l^3,s_3^3s_4^3s_5^3\right)$, see the left diagram in Figure~\ref{fig:planar}. Let $\gamma_{t_4}$ be the arc in $\omega$ which satisfies $l^4 := l(\gamma_{t_4}) = 1$ and $r(\gamma_{t_4}) = 3$, and denote $\wrd(\gamma_{t_4}) = \left(l^4,s_2^4 s_3^4\right)$. If $s_2^4 = \text{o}$, then $\gamma_{t_4}$ and $\gamma_{t_2}$ satisfy condition (2) in Lemma~\ref{lem:crossing_exists} for $i = 2$ and $j = 3$. But then any representatives of $\gamma_{t_2}$ and $\gamma_{t_4}$ will have a nontrivial crossing, a contradiction. Similarly, if $s_2^4 = \text{u}$, then $\gamma_{t_3}$ and $\gamma_{t_4}$ satisfy condition (2) in Lemma~\ref{lem:crossing_exists} for $i = 2$ and $j = 3$, again leading to a contradiction.

    Next, suppose $\omega$ contains an arc $\gamma_{t_3}$ with $\wrd(\gamma_{t_3}) = (2,\text{ooe})$, see the right diagram of Figure~\ref{fig:planar}. The argument is similar to the previous case, although more arcs must be added to deduce a contradiction. We opt to describe each arc in words, rather than using the data coming from $\wrd(-)$. First note that the arc $\gamma_{t_4}$ with left endpoint 1 and right endpoint 3 must cross under the node 2 to avoid a nontrivial crossing with either $\gamma_{t_1}$ or $\gamma_{t_2}$. Similarly, the arc $\gamma_{t_5}$ with left endpoint 1 and right endpoint 4 must pass under the nodes 2 and 3 for the same reason. It then follows that the arc $\gamma_{t_6}$ with left endpoint 0 and right endpoint 4 must pass below the nodes 1, 2, and 3 to avoid having nontrivial crossings with the previous arcs. Similarly, the arc $\gamma_{t_7}$ with left endpoint 1 and right endpoint 5 will need to pass above the nodes 2, 3, and 4 to avoid the existence of nontrivial crossings. But then the arc $\gamma_{t_8}$ with left endpoint 0 and right endpoint 3 must have a nontrivial crossing with one of the previous arcs. See Figure~\ref{fig:planar}.

    The last possibility is that $\omega$ contains an arc $\gamma_{t_3}$ with $\wrd\left(\gamma_{t_3}\right) = \left(2,\mathrm{oue}\right)$. One can then obtain a contradiction using an argument similar to the above paragraph. Precisely, the named arcs can be considered in the same order, but with $\gamma_{t_5}$, $\gamma_{t_6}$, and $\gamma_{t_7}$ mirrored along the horizontal axis.
\end{proof}

For the convenience of the reader, we include a proof of the following well-known result. Our proof is based on that of \cite[Corollary~10.21]{graph_book}, which establishes a weaker upper bound without the bipartite assumption.

\begin{lemma}\label{lem:planar_bipartite}
    Let $G$ be a simple planar bipartite graph on $(n+1) > 2$ nodes. Then $G$ has at most $2n-2$ edges.
\end{lemma}

\begin{proof}
    It suffices to consider only the case where $G$ is connected. Let $v$, $e$, and $f$ denote the number of nodes, edges, and faces of $G$, respectively. Recall that each edge of $G$ is incident to precisely two (not necessarily distinct) faces of $G$. Moreover, since $v \geq 3$ and $G$ is simple and bipartite, each face of $G$ is incident to at least four (not necessary distinct) edges. (These edges describe a closed walk in $G$, and any closed walk has even length since $G$ is bipartite.) Together, these observations yield $4f \leq 2e$. Combining this with Euler's formula ($v - e + f = 2$), we conclude that $v - e/2 \leq 2$, or equivalently that $e \leq 2n-2$.
\end{proof}

As a consequence of Proposition~\ref{prop:planar} and Lemma~\ref{lem:planar_bipartite}, we obtain the following.

\begin{corollary}\label{cor:planar}
    The maximum length of a clockwise-ordered arc diagram on $n+1$ nodes is precisely $\max\{1,2n-2\}$.
\end{corollary}

\begin{figure}
\begin{tikzpicture}[scale=0.8]

\draw[blue,very thick,smooth]
    (-2,0) [out = 45,in = 180] to (3,2) [out = 0,in = 120] to (8,0);
\draw[orange,very thick,smooth,dashed]
    (0,0) [out = 45,in = 180] to (4,1.5) [out = 0,in = 135] to (8,0);
\draw[blue,very thick,smooth]
    (2,0) [out = 45,in = 180] to (5,1) [out = 0,in = 150] to (8,0);
\draw[orange,very thick,smooth,dashed]
    (4,0) [out = 45,in = 180] to (6,0.5) [out = 0,in = 180] to (8,0);

\draw[orange,very thick,smooth,dashed]
    (-2,0) [out = -45,in = 180] to (2,-1.75) [out = 0,in = -120] to (6,0);
\draw[blue,very thick,smooth]
    (0,0) [out = -45,in = 180] to (3,-1.25) [out = 0,in = -135] to (6,0);
\draw[orange,very thick,smooth,dashed]
    (2,0) [out = -45,in = 180] to (4,-0.75) [out = 0,in = -150] to (6,0);
\draw[blue,very thick,smooth]
    (4,0)  to (6,0);

\node [orange] at (7,0) {$\gamma_1$};
\node [blue] at (3,1) {$\gamma_2$};
\node [orange] at (1,1.1) {$\gamma_3$};
\node [blue] at (-1,1.25) {$\gamma_4$};
\node [blue] at (4.5,-0.25) {$\gamma_8$};
\node [orange] at (3,-0.3) {$\gamma_7$};
\node [blue] at (1,-0.45) {$\gamma_6$};
\node [orange] at (-1,-0.5) {$\gamma_5$};

\filldraw 
 (0,0) circle (2.4pt)
 (-2,0) circle (2.4pt)
 (2,0) circle (2.4pt)
 (4,0) circle (2.4pt) 
 (6,0) circle (2.4pt) 
 (8,0) circle (2.4pt) 
 ;
 
\end{tikzpicture}
\caption{A realization of the complete bipartite graph $K_{4,2}$ as a clockwise-ordered arc diagram on 6 nodes.}\label{fig:maximal}
\end{figure}
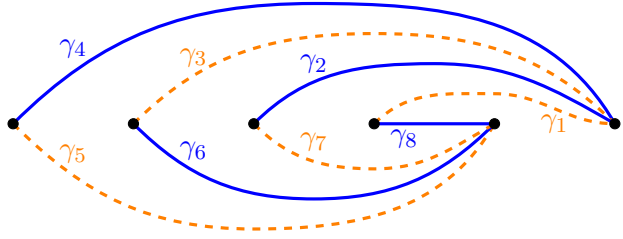

\begin{proof}
    If $n = 1$, then $\arc(n)$ contains only a single arc and there is nothing to show. Thus suppose $n > 1$. Then the fact that $2n - 2$ is an upper bound on the length of a clockwise-ordered arc diagram is an immediate consequence of Proposition~\ref{prop:planar} and Lemma~\ref{lem:planar_bipartite}.

    To see that this bound is realized, we will show that the complete bipartite graph $K_{n-1,2}$, which has precisely $2n-2$ edges, can be realized as the graph of some clockwise-ordered arc diagram. Indeed, for $1 \leq t \leq n-1$, let $\text{o}^t$ be the word consisting of $t$ copies of the letter $\text{o}$, and let $\gamma_t$ denote an arc with
    $$\wrd(\gamma_t) = (n-t-1,\text{o}^te).$$
    That is, the left endpoint of $\gamma_t$ is $n-t-1$, the right endpoint is $n$, and $\gamma_t$ passes over every node in its arrow support. Likewise for $n \leq t \leq 2n-2$, let $\text{u}^{t-n}$ be the (possibly empty) word consisting of $t$ copies of the letter $\text{u}$, and let $\gamma_t$ denote an arc with
    $$\wrd(\gamma_t) = (t-n, \text{u}^{t-n}e).$$
    That is, the left endpoint of $\gamma_t$ is $n-t$, the right endpoint is $n-1$, and $\gamma_t$ passes under every node in its arrow support. It is straighforward to show that $(\gamma_t)_{t = 2n-2}^1 \in \arc_{cw}(n,2n-2)$. See Figure~\ref{fig:maximal} for an example.
\end{proof}

We now prove the main result of this section.

\begin{theorem}\label{thm:weak}\
    Let $S \subseteq [2,n-1]_\mathbb{Z}$ and $k \in \mathbb{N}$. Then $\sigma$ induces a bijection $\weak(\Pi(A_n,S),k) \rightarrow \arc_{cw}(n,k)$; that is, the association $(X_k,\ldots,X_1) \mapsto (\sigma(X_k),\ldots,\sigma(X_1))$ is a bijection between $\weak(\Pi(A_n,S),k)$ and $\arc_{cw}(n,k)$.
\end{theorem}

\begin{proof}
    By Corollary~\ref{cor:same_weak}, it suffices to prove the result for $S = \emptyset$. Let $(X_k,\ldots,X_1) \in \weak(\Pi(A_n),k)$. Choose representatives of the arcs $\sigma(X_1),\ldots,\sigma(X_k)$ so that each pair of arcs intersect a minimal number of times. For $1 \leq i < j \leq k$, Corollary~\ref{cor:CB} and the assumption that $\Ext^1_{\Pi(A_n)}(X_i,X_j) = 0$ imply that $\sigma(X_i)$ and $\sigma(X_j)$ either have a shared endpoint or do not intersect. Moreover, since $\Hom_{\Pi(A_n)}(X_i,X_j) = 0$, Proposition~\ref{prop:hom_space} and Corollary~\ref{cor:hom} imply that $\sigma(X_i)$ and $\sigma(X_j)$ have at most one shared endpoint, and that such a shared endpoint must be directed from $\sigma(X_j)$ to $\sigma(X_i)$. We conclude that $(\sigma(X_k),\ldots,\sigma(X_1)) \in \arc_{cw}(n,k)$.
    
    Conversely, let $(\gamma_k,\ldots,\gamma_1) \in \arc_{cw}(n,k)$. The fact that $(\sigma^{-1}(\gamma_k),\ldots,\sigma^{-1}(\gamma_1)) \in \weak(\Pi(A_n),k)$ then follows immediately from Proposition~\ref{prop:hom_space}, Corollary~\ref{cor:CB}, and Corollary~\ref{cor:bricks_K}.
\end{proof}

\begin{remark}
    Since every indecomposable module over $\Pi(A_n,[2,n-1]_\mathbb{Z})$ is a brick, we note that Theorem~\ref{thm:weak} gives a classification of all stratifying systems over this algebra.
\end{remark}

Combining Theorem~\ref{thm:weak} with Corollary~\ref{cor:planar}, we obtain the following.

\begin{corollary}\label{cor:planar2}
    Let $S \subseteq [2,n-1]_\mathbb{Z}$. Then the maximum length of a weak exceptional sequence over $\Pi(A_nS)$ is $\max\{1,2n-2\}$.
\end{corollary}

We conclude this section with an example. Note in particular that (2) is an example of a weak exceptional sequence which cannot be extended into a weak exceptional sequence of length $2n-2$.

\begin{figure}
\begin{tikzpicture}[scale=0.8]

\draw[blue,very thick,smooth]
    (-2,0) [out = 45,in = 180] to (1,1.5) [out = 0,in = 135] to (4,0);

\draw[purple,very thick, smooth,dotted]
    (-2,0) [out = 30,in = 180] to (0,1) [out = 0,in = 150] to (2,0);

\draw[orange,very thick,smooth,dashed]
    (0,0) to (2,0);

\draw[purple,very thick,smooth,dotted]
    (0,0) [out = -45,in = 180] to (2,-1) [out = 0,in = -135] to (4,0);

\node [blue] at (2.5,0.75) {$\gamma_3$};
\node [purple] at (0,0.75) {$\gamma_4$};
\node [orange] at (1,-0.25) {$\gamma_1$};
\node [purple] at (3,-0.5) {$\gamma_2$};

\filldraw  (0,0) circle (2.4pt)
 (-2,0) circle (2.4pt)
 (2,0) circle (2.4pt)
 (4,0) circle (2.4pt) 
 ;
\end{tikzpicture}
\caption{A clockwise-ordered arc diagram on 4 nodes. The arc $\gamma_1$ is dashed orange, the arcs $\gamma_2$ and $\gamma_4$ are dotted purple, and the arc $\gamma_3$ is solid blue.}\label{fig:weak_ex}
\end{figure}
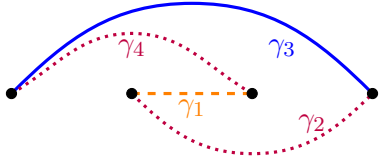

\begin{example}\label{ex:weak}\
\begin{enumerate}
    \item Consider the clockwise-ordered arc diagram $(\gamma_4,\gamma_3,\gamma_2,\gamma_1)$ shown in Figure~\ref{fig:weak_ex}. The corresponding weak exceptional sequence in either $\mods\Pi(A_n,[2,n-1]_\mathbb{Z})$ or $\mods \Pi(A_n)$ is $$\left({\footnotesize \begin{matrix}2\\1\end{matrix}}, \ {\scriptsize \begin{matrix}3\\2\\1\end{matrix}}, \ {\footnotesize\begin{matrix}2\\3\end{matrix}}, \ 2\right).$$
    
    \item Consider the clockwise-ordered arc diagram in Figure~\ref{fig:maximal}. By deleting the rightmost node (number 5) and all incident arcs, we obtain a clockwise-ordered arc diagram $(\gamma_8,\gamma_7,\gamma_6,\gamma_5)$. The corresponding weak exceptional sequence over either $\Pi(A_4,[2,n-1]_\mathbb{Z})$ or $\Pi(A_4)$ is $$\left(4,\ {\footnotesize\begin{matrix}3\\4\end{matrix}},\ {\scriptsize \begin{matrix}2\\3\\4\end{matrix}},\ {\tiny\begin{matrix}1\\2\\3\\4\end{matrix}}\right).$$
    Note that while this clockwise-ordered arc diagram consists of $4 = n < 2n-2$ arcs, there are no additional arcs that can be added anywhere in the sequence without violating the clockwise-ordered property.
    \item Consider the clockwise-ordered arc diagram $(\gamma_3,\gamma_2,\gamma_1)$ shown in Figure~\ref{fig:clockwise}. The corresponding weak exceptional sequence over either $\Pi(A_5,[2,n-1]_\mathbb{Z})$ or $\Pi(A_5)$ is 
    $$\left({\scriptsize \begin{matrix} \ 2\\1 \ 3\\\phantom{44}4\end{matrix}}, \ {\scriptsize\begin{matrix}\phantom{44} 3 \ 5\\ 2 \ 4\\1\ \phantom{44}\end{matrix}},\ {\scriptsize \begin{matrix}3\\2\\1\end{matrix}}\right).$$
\end{enumerate}
It is shown in the concurrent work \cite{hanson} that the sequences in (2) and (3) above are examples of ``$\tau$-exceptional sequences'', while the one in (1) is not.
\end{example}

\section{Exceptional sequences over hereditary algebras of type A}\label{sec:hereditary}

In this section, we briefly explain how the results of this paper relate to the model of exceptional sequences over hereditary algebras of type A constructed in \cite{GIMO}.

Let $\varepsilon: \{1,\ldots,n-1\} \rightarrow \{\textnormal{u},\textnormal{o}\}$. The function $\varepsilon$ defines a quiver $A_n^\varepsilon$, which, for any $S \subseteq [2,n-1]_\mathbb{Z}$, we can see as a quotient of $\Pi(A_n,S)$ as follows. Let $E(\varepsilon) = \{a_i \mid \varepsilon(i) = \textnormal{u}\} \cup \{a_i^* \mid \varepsilon(i) = \textnormal{o}\}$, and let $KA_n^\varepsilon = K\overline{A_n}/(E(\varepsilon))$. Alternatively, $KA_n^\varepsilon$ is the path algebra of a quiver of type $A_n$ with an arrow $a_i: i \rightarrow i+1$ whenever $\varepsilon(i) = \textnormal{o}$ and an arrow $a_i^*: i+1 \rightarrow i$ whenever $\varepsilon(i) = \textnormal{u}$.

Let $\gamma \in \arc(n)$ and write $\wrd(\gamma) = \left(l(\gamma),s_{l(\gamma)+1}\cdots s_{r(\gamma)}\right)$. We say that $\gamma$ is \emph{$\varepsilon$-admissible} if $s_i = \varepsilon(i)$ for all $l(\gamma) < i < r(\gamma)$. It is well-known that $\sigma$ induces a bijection from $\brick(KA_n^\varepsilon)$ to the set of (combinatorial equivalence classes of) $\varepsilon$-admissible arcs.

\begin{remark}
    One may also see a choice of $\varepsilon$ as a choice of ``Coxeter element'' $c$ in the symmetric group $\mathfrak{S}_{n+1}$. The $\varepsilon$-admissible arcs can then be seen as those which are ``$c$-sortable''. See e.g. \cite[Section~4]{mizuno2} or \cite{BR} for further discussion.
\end{remark}

In \cite{GIMO}, Garver, Igusa, Mathern, and Ostroff classify the exceptional sequences (which in this case correspond precisely to the weak exceptional sequences) over $KA_n^\varepsilon$ using ``labeled strand diagrams''. Their definition includes the assumption that the arcs underlying any strand diagram form a forest. We remove this assumption and restate the remainder of their definition in our terminology as follows.

\begin{definition}\textnormal{(c.f. \cite[Definiton~14]{GIMO})}
    Let $\Delta = (\gamma_k,\ldots,\gamma_1)$ be a sequence of $\varepsilon$-admissible arcs. We say that $\Delta$ is an \emph{$\varepsilon$-admissible labeled strand diagram} if 
    the following hold for all $1 \leq i < j \leq k$.
        \begin{enumerate}
            \item There exist representatives of $\gamma_i$ and $\gamma_j$ which do not have a nontrivial crossing.
            \item If $\gamma_i$ and $\gamma_j$ have a shared endpoint then $\gamma_j$ is clockwise from $\gamma_i$.
            \item Suppose that $\gamma_i$ and $\gamma_j$ have a contested endpoint $x$. If $\varepsilon(x) = \textnormal{o}$ then $x = r(\gamma_i)$. Similarly, if $\varepsilon(x) = \textnormal{u}$ then $x = r(\gamma_j)$.
    \end{enumerate}
    We denote by $\arc_\varepsilon(n,k)$ the set of $\varepsilon$-admissible labeled strand diagrams and $\arc_\varepsilon(n) = \bigcup_{k \in \mathbb{N}} \varepsilon(n,k).$
\end{definition}

We now state and prove the classification given in \cite{GIMO}, but again without the assumption that the arcs in question form a forest. Recall once again that $\weak(KA_n^\varepsilon)$ can be considered as the set of exceptional sequences of $KA_n^\varepsilon$.

\begin{theorem}\textnormal{(c.f. \cite[Theorem~16]{GIMO})}\label{thm:hereditary}
 The map $\sigma$ induces a bijection $\weak(KA_n^\varepsilon,k) \rightarrow \arc_\varepsilon(n,k)$.
\end{theorem}

\begin{proof}
    It suffices to prove the result for $k = 2$.
    Denote by $F: \mods KA_n^\varepsilon \rightarrow \mods \Pi(A_n,[2,n-1]_\mathbb{Z})$ the inclusion functor induced by the quotient map $\Pi(A_n,[2,n-1]_\mathbb{Z}) \rightarrow KA_n^\varepsilon$. Note that, as with the functors $F_T^S$, this $F$ is a fully faithful functor which induces injective linear maps on $\Ext^1$. Moreover, we can see the application of this functor as turning representations of $A_n^\varepsilon$ into representations of $\overline{A_n}$ by adding the 0 map on all new arrows. 
    
    Suppose first that $\sigma(X)$ and $\sigma(Y)$ have a contested endpoint $x$. Then $\Hom_{KA_n^\varepsilon}(X,Y) =0$ and $\dim\Ext^1_{\Pi(A_n,[2,n-1]_\mathbb{Z})}(FX,FY) = 1$ by Proposition~\ref{prop:hom_space} and Theorem~\ref{thm:ext_space}. Now let $\eta_1(x)$ be the short exact sequence from Equation~\eqref{extension1}, which spans $\Ext^1_{\Pi(A_n,[2,n-1]_\mathbb{Z})}(FX,FY)$ by Proposition~\ref{prop:extension1}. The middle term of this short exact sequence is a module over $KQ_n^\varepsilon$ if and only if either (a) $\varepsilon(x) = \textnormal{u}$ and $x = r(\sigma(X))$ or (b) $\varepsilon(x) = \textnormal{o}$ and $x = r(\sigma(Y))$. This shows that $(Y,X) \in \weak(KA_n^\varepsilon,2)$ if and only if $(\sigma(Y),\sigma(X)) \in \arc_\varepsilon(n,2)$ in this case.

    Now suppose that $X$ and $Y$ do not have a contested endpoint. Then $\Ext^1_{\Pi(A_n,[2,n-1]_\mathbb{Z})}(FX,FY)$ is spanned by the short exact sequences of the form $\eta_2(x)$ in Equation~\eqref{extension2}. Moreover, since $FX$ and $FY$ correspond to $\varepsilon$-admissible arcs, the middle terms of these short exact sequences will, as well. Thus $(FY,FX) \in \weak(\Pi(A_n,[2,n-1]_\mathbb{Z})$ if and only if $(Y,X) \in \weak(KA_n^\varepsilon)$. Lemma~\ref{lem:same_weak} thus implies the result in this case.
\end{proof}

\begin{remark}
    Our proof of Theorem~\ref{thm:hereditary} highlights that it is not necessary to assume that the underlying arcs form a forest in \cite[Definition~14]{GIMO}. Indeed, it is shown in the proof of \cite[Theorem~12]{GIMO} that the arcs corresponding to any exceptional sequence do form a forest, but the forest condition is not used in showing that the modules corresponding to any strand diagram form an exceptional sequence.
\end{remark}

We conclude by highlighting two critical differences between weak exceptional sequences over $\Pi(A_n,S)$ and exceptional sequences over $\Pi(A_n^\varepsilon)$.

\begin{example}\
    \begin{enumerate}
        \item Suppose $(X_k,\ldots,X_1)$ is an exceptional sequence over $KA_n^\varepsilon$. Then it may not be the case that $(FX_k,\ldots,FX_1)$ is an exceptional sequence over $\Pi(A_n,S)$. Indeed, let $\gamma_1$ and $\gamma_2$ be as in Figure~\ref{fig:crossings}, and suppose that $\varepsilon(5) = \textnormal{o}$. Then $(\sigma^{-1}(\gamma_2),\sigma^{-1}(\gamma_1))$ is an exceptional sequence over $KA_n^\varepsilon$, but is not a weak exceptional sequence over $\Pi(A_n,S)$ (for any $S$).
        \item Example~\ref{ex:weak}(1) gives an example of a weak exceptional sequence over $\Pi(A_n,S)$ (for any $S$) whose corresponding arcs do not form a forest. The argument used to prove the forest result over $KA_n^\varepsilon$ in \cite{GIMO} goes wrong in this example as follows. Traversing the arcs in right-to-left order, the cycle $(\gamma_2,\gamma_1,\gamma_4,\gamma_3)$ is oriented counterclockwise around the region it contains. Identifying indices mod 4, the argument in \cite[Theorem~12]{GIMO} then deduces that $\gamma_{i+1}$ is clockwise of $\gamma_i$ for all $i$. The reason this fails for $i = 4$ in our example is that $\gamma_2$ is allowed to cross under the node 2, even though $\gamma_3$ has already passed above this node.
    \end{enumerate}
\end{example}


\bibliographystyle{amsalpha}
\bibliography{biblio.bib}

\end{document}